\documentclass[11pt]{article}
\usepackage{setspace}
\usepackage{amsmath} 
\usepackage{amssymb}
\usepackage{latexsym}
\usepackage{xcolor}
\usepackage{fancyhdr}
\allowdisplaybreaks 

 \usepackage{comment}

\def\XXint#1#2#3{{\setbox0=\hbox{$#1{#2#3}{\int}$} 
\vcenter{\hbox{$#2#3$}}\kern-.5\wd0}}   

 \numberwithin{equation}{section}
\newtheorem{theorem}[equation]{Theorem}
\newtheorem{proposition}[equation]{Proposition}

\newtheorem{lemma}[equation]{Lemma}

\title{On the tangential gradient of the kernel of the double layer potential}

\author{  
Massimo Lanza de Cristoforis
\\
Dipartimento di Matematica `Tullio Levi-Civita', 
\\
Universit\`a degli Studi di Padova, 
\\
Via Trieste 63, Padova 35121, 
Italy. 
\\
E-mail: mldc@math.unipd.it   }

 \date{\ }

\begin{document} 
\maketitle

\noindent
{\bf Abstract:}  In this paper we consider an elliptic operator with constant coefficients and 
we estimate the maximal function of the tangential gradient of the kernel of the double layer potential with respect to its first variable. As a consequence,  we deduce the validity of a continuity property of the double layer potential in H\"{o}lder spaces on the boundary that extends previous results for the Laplace operator and for the Helmholtz operator.  

 \vspace{\baselineskip}

\par


\section{Introduction}
 In this paper, we consider the kernel of the double layer potential associated to the fundamental solution of a second order differential operator with constant coefficients, that we now introduce. Unless otherwise specified,  we assume  throughout the paper that
\[
n\in {\mathbb{N}}\setminus\{0,1\}\,,
\]
where ${\mathbb{N}}$ denotes the set of natural numbers including $0$.  Let $\Omega$ be a bounded open subset of ${\mathbb{R}}^{n}$ of class $C^{1}$.   We employ the same notation of reference \cite{DoLa17} with Dondi that we now introduce. 

Let $\nu_\Omega$ or simply 
  $\nu \equiv (\nu_{l})_{l=1,\dots,n}$ denote the external unit normal to $\partial\Omega$.  Let $N_{2}$ denote the number of multi-indexes $\gamma\in {\mathbb{N}}^{n}$ with $|\gamma|\leq 2$. For each 
\begin{equation}
\label{introd0}
{\mathbf{a}}\equiv (a_{\gamma})_{|\gamma|\leq 2}\in {\mathbb{C}}^{N_{2}}\,, 
\end{equation}
we set 
\[
a^{(2)}\equiv (a_{lj} )_{l,j=1,\dots,n}\qquad
a^{(1)}\equiv (a_{j})_{j=1,\dots,n}\qquad
a\equiv a_{0}\,.
\]
with $a_{lj} \equiv 2^{-1}a_{e_{l}+e_{j}}$ for $j\neq l$, $a_{jj} \equiv
 a_{e_{j}+e_{j}}$,
and $a_{j}\equiv a_{e_{j}}$, where $\{e_{j}:\,j=1,\dots,n\}$  is the canonical basis of ${\mathbb{R}}^{n}$. We note that the matrix $a^{(2)}$ is symmetric. 
Then we assume that 
  ${\mathbf{a}}\in  {\mathbb{C}}^{N_{2}}$ satisfies the following ellipticity assumption
\begin{equation}
\label{ellip}
\inf_{
\xi\in {\mathbb{R}}^{n}, |\xi|=1
}{\mathrm{Re}}\,\left\{
 \sum_{|\gamma|=2}a_{\gamma}\xi^{\gamma}\right\} >0\,,
\end{equation}
and we consider  the case in which
\begin{equation}
\label{symr}
a_{lj} \in {\mathbb{R}}\qquad\forall  l,j=1,\dots,n\,.
\end{equation}
Then we introduce the operators
\begin{eqnarray*}
P[{\mathbf{a}},D]u&\equiv&\sum_{l,j=1}^{n}\partial_{x_{l}}(a_{lj}\partial_{x_{j}}u)
+
\sum_{l=1}^{n}a_{l}\partial_{x_{l}}u+au\,,
\\
B_{\Omega}^{*}v&\equiv&\sum_{l,j=1}^{n} \overline{a}_{jl}\nu_{l}\partial_{x_{j}}v
-\sum_{l=1}^{n}\nu_{l}\overline{a}_{l}v\,,
\end{eqnarray*}
for all $u,v\in C^{2}(\overline{\Omega})$, and a fundamental solution $S_{{\mathbf{a}} }$ of $P[{\mathbf{a}},D]$, and the  boundary integral operator corresponding to the  double layer potential 
\begin{eqnarray}
\label{introd3}
\lefteqn{
W_\Omega[{\mathbf{a}},S_{{\mathbf{a}}}   ,\mu](x) \equiv 
\int_{\partial\Omega}\mu (y)\overline{B^{*}_{\Omega,y}}\left(S_{{\mathbf{a}}}(x-y)\right)
\,d\sigma_{y}
}
\\  \nonumber
&&
\qquad
=-\int_{\partial\Omega}\mu(y)\sum_{l,j=1}^{n} a_{jl}\nu_{l}(y)\frac{\partial S_{ {\mathbf{a}} } }{\partial x_{j}}(x-y)\,d\sigma_{y}
\\  \nonumber
&&
\qquad\quad
-\int_{\partial\Omega}\mu(y)\sum_{l=1}^{n}\nu_{l}(y)a_{l}
S_{ {\mathbf{a}} }(x-y)\,d\sigma_{y} \qquad\forall x\in \partial\Omega\,,
\end{eqnarray}
where the density or moment $\mu$ is a function  from $\partial\Omega$ to ${\mathbb{C}}$. Here the subscript $y$ of $\overline{B^{*}_{\Omega,y}}$ means that we are taking $y$ as variable of the differential operator $\overline{B^{*}_{\Omega,y}}$. The kernel of the double layer potential on $\partial\Omega$ is the following
\begin{equation}\label{eq:tgdlgen}
\overline{B^{*}_{\Omega,y}}\left(S_{{\mathbf{a}}}(x-y)\right)\equiv - \sum_{l,j=1}^{n} a_{jl}\nu_{l}(y)\frac{\partial S_{ {\mathbf{a}} } }{\partial x_{j}}(x-y) 
 - \sum_{l=1}^{n}\nu_{l}(y)a_{l}
S_{ {\mathbf{a}} }(x-y) 
\end{equation}
for all $(x,y)\in (\partial\Omega)^2\setminus {\mathbb{D}}_{\partial\Omega}$,
where
$
{\mathbb{D}}_{\partial\Omega}\equiv\{(x,y)\in (\partial\Omega)^2:\,x=y\} 
$ denotes the diagonal of $ (\partial\Omega)^2$. 

The role of the double layer potential in the solution of boundary value problems for the operator $P[{\mathbf{a}},D]$ is well known (cf.~\textit{e.g.}, 
G\"{u}nter~\cite{Gu67}, Kupra\-dze,  Gegelia,  Basheleishvili and 
 Burchuladze~\cite{KuGeBaBu79}, Mikhlin \cite{Mik70},  Mikhlin and  Pr\"{o}ssdorf \cite{MikPr86}).	 
  In order to prove that under suitable assumptions  $W_\Omega[{\mathbf{a}},S_{{\mathbf{a}}}   ,\cdot]$ maps a H\"{o}lder space in $\partial\Omega$  to a Schauder space of differentiable functions on $\partial\Omega$, one needs to estimate the maximal function of the tangential gradient with respect to the first variable of the kernel of the double layer potential. More precisely, one needs to prove that
 the maximal function of the tangential gradient of the kernel of the double layer potential with respect to its first variable is bounded,  \textit{i.e.}, 	that
\begin{equation}\label{thm:dllreggenn1}
\sup_{x\in \partial\Omega}\sup_{\rho\in ]0,+\infty[}
\left\vert
\int_{(\partial\Omega)\setminus {\mathbb{B}}_{n}(x,\rho)} {\mathrm{grad}}_{\partial\Omega,x}\overline{B^{*}_{\Omega,y}}\left(S_{{\mathbf{a}}}(x-y)\right)\,d\sigma_y
\right\vert <+\infty
\end{equation}
where 
\[
{\mathbb{B}}_{n}(x,r)\equiv\{y\in {\mathbb{R}}^n:\,|y-x|<r\}
\]
for all $x\in {\mathbb{R}}^n$ and $r\in]0,+\infty[$.  See for example Colton and Kress \cite[(2.26)]{CoKr83}, Kirsch and Hettlich \cite[(3.25e)]{KiHe15} for the Helmoltz operator and for applications to the equations of the electromagnetism and Hsiao and Wendland \cite[Remark 1.2.1]{HsWe08} for the Laplace operator in case $\Omega$ is of class $C^2$. The aim of this paper is to show that if $\Omega$ is of class $C^{1,1}$, then such a maximal function is bounded for a general second order differential  operator as $P[{\mathbf{a}},D]$ 
(see Theorem \ref{thm:mftgdlgx}). Then as a consequence of inequality (\ref{thm:dllreggenn1}), the following result holds true (see  \cite{La22d}).
\begin{theorem}\label{thm:dllreggen}
  Let  $\beta\in]0,1]$.  Let $\Omega$ be a bounded open subset of ${\mathbb{R}}^{n}$ of class  $C^{1,1}$.  
    Let ${\mathbf{a}}$ be as in (\ref{introd0}), (\ref{ellip}), (\ref{symr}).  Let $S_{ {\mathbf{a}} }$ be a fundamental solution of $P[{\mathbf{a}},D]$. Then the following statements hold.
  \begin{enumerate}
\item[(i)] If $\beta<1$, then the  operator $W_\Omega[{\mathbf{a}},S_{{\mathbf{a}}}   ,\cdot]$ from 
$C^{0,\beta}(\partial\Omega)$ to $C^{1, \beta }(\partial\Omega)$ defined by 
 (\ref{introd3}) for all $\mu\in C^{0,\beta}(\partial\Omega)$ is linear and continuous.
 \item[(ii)] If $\beta=1$, then the  operator $W_\Omega[{\mathbf{a}},S_{{\mathbf{a}}}   ,\cdot]$ from 
$C^{0,1}(\partial\Omega)$ to $C^{1, \omega_1(\cdot) }(\partial\Omega)$ defined by 
 (\ref{introd3})  for all $\mu\in C^{0,1}(\partial\Omega)$ is linear and continuous.
\end{enumerate}
\end{theorem}
Here $C^{1,\omega_{1}}(\partial\Omega)$ denotes the generalized Schauder space
  of functions with $1$-st order tangential  derivatives that satisfy a generalized $\omega_{1}$-H\"{o}lder condition with
\[
\omega_{1}(r)
\equiv
\left\{
\begin{array}{ll}
0 &r=0\,,
\\
r |\ln r | &r\in]0,r_{1}]\,,
\\
r_{1} |\ln r_{1} | & r\in ]r_{1},+\infty[\,,
\end{array}
\right.
\]
where
$
r_{1}\equiv e^{-1}
$. For the classical definition of the generalized H\"{o}lder or Schauder spaces we refer the reader to
  \cite[p.~75, \S 2]{DoLa17} and to reference \cite[\S 2.6, \S 2.11]{DaLaMu21} with Dalla Riva and Musolino. Thus Theorem \ref{thm:dllreggen} is a generalization of   classical results  for the Laplace and Helmoltz operator in case $\Omega$ is of class $C^2$ and $\beta\in]0,1[$ (see Colton and Kress \cite[Thm. 2.22]{CoKr83}, 
Hsiao and Wendland \cite[Remark 1.2.1]{HsWe08}).

\section{Technical preliminaries on the differential operator}\label{sec:techdiop}
  Let $M_n({\mathbb{R}})$ denote the set of $n\times n$ matrices with real entries. 
     $\delta_{l,j}$ denotes the Kronecker\index{Kronecker symbol}  symbol. Namely,  $\delta_{l,j}=1$ if $l=j$, $\delta_{l,j}=0$ if $l\neq j$, with $l,j\in {\mathbb{N}}$. $|A|$ denotes the operator norm of a matrix $A$, 
       $A^{t}$ denotes the transpose matrix of $A$.	 Let $O_{n}({\mathbb{R}})$ denote the set of $n\times n$ orthogonal matrices with real entries.  
 Let $\Omega$ be an open subset of ${\mathbb{R}}^n$. Let $s\in {\mathbb{N}}\setminus\{0\}$, $f\in \left(C^{1}(\Omega)\right)^{s} $. Then   $Df$ denotes the Jacobian matrix of $f$. We say that $\Omega$ is of class $C^0$ provided that for each point $p\in\partial\Omega$ there exist
$R_p\in O_{n}({\mathbb{R}})$, $r$, $\delta\in ]0,+\infty[$ such that the intersection
 \[
 R_p(\Omega-p) \cap ({\mathbb{B}}_{n-1}(0,r)\times]-\delta,\delta[    )
 \]
 is the strict hypograph of a continuous function $\gamma_p$ from ${\mathbb{B}}_{n-1}(0,r)$ to $]-\delta ,\delta [$ which vanishes at $0$ and such that $|\gamma_p(\eta)|<\delta/2$ for all $\eta\in {\mathbb{B}}_{n-1}(0,r)$,  i.e., provided that there exists $\gamma_p\in C^{0}({\mathbb{B}}_{n-1}(0,r), ]-\delta ,\delta [)$ such that 
 \begin{eqnarray}
\label{prelim.cocylind1}
\lefteqn{R_p(\Omega-p )\cap ({\mathbb{B}}_{n-1}(0,r)\times ]-\delta,\delta[) 
}
\\ \nonumber
&&\qquad 
=\left\{
(\eta,y)\in 
{\mathbb{B}}_{n-1}(0,r)\times ]-\delta,\delta[:\, y<\gamma_p(\eta)
\right\}
\equiv{\mathrm{hypograph}}_{s}(\gamma_p) 
\,, 
\\ \nonumber
&&\qquad 
|\gamma_p(\eta)|<\delta/2\qquad\forall \eta\in {\mathbb{B}}_{n-1}(0,r)\,,\qquad
\gamma_p(0)=0\,.
\end{eqnarray}
Here the subscript `$s$' of  `$\,{\mathrm{hypograph}}_{s}$' stands for `strict'. Then we say that the set
 \[
 C(p,R_p,r,\delta)\equiv p+ R_p^{t}({\mathbb{B}}_{n-1}(0,r)\times]-\delta,\delta[    )
 \]
 is a coordinate cylinder  for $\Omega$ around $p$, that the function $\gamma_p$ represents $\partial\Omega$ in $C(p,R_p,r,\delta)$ and that the function $\psi_{p}$ from ${\mathbb{B}}_{n-1}(0,r)$ to ${\mathbb{R}}^{n}$ defined by 
\begin{equation}
\label{prelim.cocylind3}
\psi_{p}(\eta)\equiv p+R_p^{t}\left(
\begin{array}{c}
\eta
\\
\gamma_p(\eta)
\end{array}
\right)\qquad\forall \eta\in {\mathbb{B}}_{n-1}(0,r)
\end{equation} 
 is the parametrization of $\partial\Omega$ around $p$ in the coordinate cylinder $C(p,R_p,r,\delta)$.  Then we say that $\Omega$ is of class $C^{m}$ or of class $C^{m,\alpha}$
for some $m\in{\mathbb{N}}$, $\alpha\in ]0,1]$ provided that  $\gamma_p$ is of class $C^{m}$ or of class $C^{m,\alpha}$ for all $p\in\partial\Omega$ (as in reference \cite[\S 2.7, \S 2.13]{DaLaMu21} with Dalla Riva and Musolino). If ${\mathbb{D}}$ is a subset of ${\mathbb{R}}^n$, then  $|f:\,{\mathbb{D}}|_\alpha$ (or simply $|f|_\alpha$) denotes the $\alpha$-H\"older constant of $f$ and    ${\mathrm{Lip}}(f)$ denotes the Lipschitz constant of $f$. 
    In order to analyze the kernel of the double layer potential, we need some more information on the fundamental solution $S_{ {\mathbf{a}} } $.  To do so, we introduce the fundamental solution $S_{n}$ of the Laplace operator. Namely, we set
\[
S_{n}(x)\equiv
\left\{
\begin{array}{lll}
\frac{1}{s_{n}}\ln  |x| \qquad &   \forall x\in 
{\mathbb{R}}^{n}\setminus\{0\},\quad & {\mathrm{if}}\ n=2\,,
\\
\frac{1}{(2-n)s_{n}}|x|^{2-n}\qquad &   \forall x\in 
{\mathbb{R}}^{n}\setminus\{0\},\quad & {\mathrm{if}}\ n>2\,,
\end{array}
\right.
\]
where $s_{n}$ denotes the $(n-1)$ dimensional measure of 
$\partial{\mathbb{B}}_{n}(0,1)$ and
we follow a formulation of Dalla Riva \cite[Thm.~5.2, 5.3]{Da13} and Dalla Riva, Morais and Musolino \cite[Thm.~5.5]{DaMoMu13}, that we state  as in paper  \cite[Cor.~4.2]{DoLa17} with Dondi (see also John~\cite{Jo55}, Miranda~\cite{Mi65} for homogeneous operators, and Mitrea and Mitrea~\cite[p.~203]{MitMit13}).   
\begin{proposition}
 \label{prop:ourfs} 
Let ${\mathbf{a}}$ be as in (\ref{introd0}), (\ref{ellip}), (\ref{symr}). 
Let $S_{ {\mathbf{a}} }$ be a fundamental solution of $P[{\mathbf{a}},D]$. 
Then there exist an invertible matrix $T\in M_{n}({\mathbb{R}})$ such that
\begin{equation}
\label{prop:ourfs0}
a^{(2)}=TT^{t}\,,
\end{equation}
 a real analytic function $A_{1}$ from $\partial{\mathbb{B}}_{n}(0,1)\times{\mathbb{R}}$ to ${\mathbb{C}}$ such that
 $A_{1}(\cdot,0)$ is odd,    $b_{0}\in {\mathbb{C}}$, a real analytic function $B_{1}$ from ${\mathbb{R}}^{n}$ to ${\mathbb{C}}$ such that $B_{1}(0)=0$, and a real analytic function $C $ from ${\mathbb{R}}^{n}$ to ${\mathbb{C}}$ such that
\begin{eqnarray}
\label{prop:ourfs1}
\lefteqn{S_{ {\mathbf{a}} }(x)
= 
\frac{1}{\sqrt{\det a^{(2)} }}S_{n}(T^{-1}x)
+|x|^{3-n}A_{1}(\frac{x}{|x|},|x|)
}
\\ \nonumber
&&\qquad
 +(B_{1}(x)+b_{0}(1-\delta_{2,n}))\ln  |x|+C(x)\,,
\end{eqnarray}
for all $x\in {\mathbb{R}}^{n}\setminus\{0\}$,
 and such that both $b_{0}$ and $B_{1}$   equal zero
if $n$ is odd. Moreover, 
 \[
 \frac{1}{\sqrt{\det a^{(2)} }}S_{n}(T^{-1}x) 
 \]
is a fundamental solution for the principal part
  of $P[{\mathbf{a}},D]$.
\end{proposition}
In particular for the statement that $A_{1}(\cdot,0)$ is odd, we refer to
Dalla Riva, Morais and Musolino \cite[Thm.~5.5, (32)]{DaMoMu13}, where $A_{1}(\cdot,0)$ coincides with ${\mathbf{f}}_1({\mathbf{a}},\cdot)$ in that paper. Here we note that a function $A$ from $(\partial{\mathbb{B}}_{n}(0,1))\times{\mathbb{R}}$ to ${\mathbb{C}}$ is said to be real analytic provided that it has a real analytic extension   to an open neighbourhood of $(\partial{\mathbb{B}}_{n}(0,1))\times{\mathbb{R}}$ in 
${\mathbb{R}}^{n+1}$. Then we have the following elementary lemma    (cf.~\textit{e.g.},  \cite[\S 4]{La22d}).		 
\begin{lemma}\label{lem:anexsph}
 Let $n\in {\mathbb{N}}\setminus\{0,1\}$. A function $A$ from   $(\partial{\mathbb{B}}_{n}(0,1))\times{\mathbb{R}}$ to ${\mathbb{C}}$ is  real analytic if and only if the function $\tilde{A}$ from $({\mathbb{R}}^n\setminus\{0\}) \times{\mathbb{R}}$ defined by 
\begin{equation}\label{lem:anexsph1}
\tilde{A}(x,r)\equiv A(\frac{x}{|x|},r)\qquad\forall (x,r)\in ({\mathbb{R}}^n\setminus\{0\}) \times{\mathbb{R}}
\end{equation}
is real analytic.
 \end{lemma}
 
Then  one can prove the following formula for the gradient of the fundamental solution (see reference \cite[Lem.~4.3,  (4.8) and the following 2 lines]{DoLa17} with Dondi). Here one should remember that $A_1(\cdot,0)$ is odd and that $b_0=0$ if $n$ is odd). 
\begin{proposition}
\label{prop:grafun}
 Let ${\mathbf{a}}$ be as in (\ref{introd0}), (\ref{ellip}), (\ref{symr}). Let $T\in M_{n}({\mathbb{R}})$  be as in (\ref{prop:ourfs0}). Let $S_{ {\mathbf{a}} }$ be a fundamental solution of $P[{\mathbf{a}},D]$. Let  $B_{1}$, $C$
 be as in Proposition \ref{prop:ourfs}. 
  Then there exists a real analytic function $A_{2}$ from $\partial{\mathbb{B}}_{n}(0,1)\times{\mathbb{R}}$ to ${\mathbb{C}}^{n}$ such that
\begin{eqnarray}
\label{grafun1}
\lefteqn{
DS_{ {\mathbf{a}} }(x)=\frac{1}{ s_{n}\sqrt{\det a^{(2)} } }
|T^{-1}x|^{-n}x^{t}(a^{(2)})^{-1} 
}
\\ \nonumber
&&\qquad\qquad
+|x|^{2-n}A_{2}(\frac{x}{|x|},|x|)+DB_{1}(x)\ln |x|+DC(x)
\end{eqnarray}
for all $x\in {\mathbb{R}}^{n}\setminus\{0\}$. 
Moreover,   $A_2(\cdot,0)$ is even.
\end{proposition}
We also note that
we note that
\begin{eqnarray}\label{eq:ppgr}
\lefteqn{\frac{1}{\sqrt{\det a^{(2)} }}D(S_n\circ T^{-1} )(x)a^{(2)} 
}
\\ \nonumber
&&  \quad
=
\frac{1}{s_n\sqrt{\det a^{(2)} }}\frac{\left((T^{-1})x\right)^tT^{-1}}{|T^{-1}x|^n}a^{(2)}
=
\frac{1}{s_n\sqrt{\det a^{(2)} }}\frac{x^t(T^{-1})^tT^{-1}a^{(2)}
}{|T^{-1}x|^n} 
\\ \nonumber
&&  \quad
=
\frac{1}{s_n\sqrt{\det a^{(2)} }}\frac{x^t 
}{|T^{-1}x|^n} 
\ \ \forall x\in {\mathbb{R}}^n\setminus\{0\}\,.
\end{eqnarray}
Then we can prove the following formula 
for the kernel of the double layer potential in case $\Omega$ is at least of class $C^1$.
\begin{eqnarray}
\label{eq:boest3}
\lefteqn{
\overline{B^{*}_{\Omega,y}}\left(S_{{\mathbf{a}}}(x-y)\right)
}
\\ \nonumber
&& 
= -DS_{{\mathbf{a}}}(x-y)    a^{(2)}  \nu  (y)
- 
\nu^{t} (y)
a^{(1)}S_{{\mathbf{a}}}(x-y)
\\ \nonumber
&& 
= 
-\frac{1}{s_{n}\sqrt{\det a^{(2)} }}|T^{-1}(x-y)|^{-n} (x-y)^{t}\nu (y)
\\ \nonumber
&& \quad
-
|x-y|^{2-n}
A_{2}(\frac{x-y}{|x-y|},|x-y|)a^{(2)}  \nu  (y)
\\ \nonumber
&&  \quad
-DB_{1}(x-y)a^{(2)}  \nu  (y)\ln |x-y|
-   DC(x-y)a^{(2)}  \nu  (y)
\\ \nonumber
&&  \quad
-\nu^{t}(y) a^{(1)}S_{{\mathbf{a}}}(x-y)
\qquad\forall x,y\in\partial\Omega,  x\neq y\,.
\end{eqnarray}
(see reference \cite[(5.2) p.~86]{DoLa17} with Dondi).  By applying equality (\ref{eq:boest3}),   we can  compute a formula for the  tangential gradient with respect to the first variable of the kernel of the double layer potential. For the definition of tangential gradient $ {\mathrm{grad}}_{\partial\Omega}$, we refer   to Kirsch and Hettlich \cite[A.5]{KiHe15}, Chavel~\cite[Chap.~1]{Cha84}. 

\begin{lemma}\label{lem:tgdlgen}  
Let $n\in {\mathbb{N}}\setminus\{0,1\}$. 
 Let ${\mathbf{a}}$ be as in (\ref{introd0}), (\ref{ellip}), (\ref{symr}).  Let $S_{ {\mathbf{a}} }$ be a fundamental solution of $P[{\mathbf{a}},D]$. 
 
 Let $\alpha\in]0,1]$. Let $\Omega$ be a bounded open subset of ${\mathbb{R}}^{n}$ of class $C^{1,\alpha}$.     If  $h\in\{1,\dots,n\}$, then
\begin{eqnarray}\label{lem:tgdlgen1}
\lefteqn{
({\mathrm{grad}}_{\partial\Omega,x}\overline{B^{*}_{\Omega,y}}\left(S_{{\mathbf{a}}}(x-y)\right))_h
=\frac{\partial}{\partial x_h}\overline{B^{*}_{\Omega,y}}\left(S_{{\mathbf{a}}}(x-y)\right))
}
\\ \nonumber
&& -
\nu_h(x)\sum_{l=1}^n \nu_l(x)
\frac{\partial}{\partial x_l}\overline{B^{*}_{\Omega,y}}\left(S_{{\mathbf{a}}}(x-y)\right)) 
\\ \nonumber
&& 
=\frac{n}{s_n\sqrt{\det a^{(2)} }}
 \frac{(x-y)^t\cdot\nu(y)}{|T^{-1}(x-y)|^n}
 \\ \nonumber
&&\quad 
 \times\sum_{l=1}^n\nu_l(x)
 \biggl[
 \nu_l(x)
 \frac{
 \sum_{j,z=1}^n (T^{-1})_{jz}(x_z-y_z)(T^{-1})_{jh}
  }{|T^{-1}(x-y)|^2}  
  \\ \nonumber
&& \quad 
 -
 \nu_h(x)\frac{\sum_{j,z=1}^n (T^{-1})_{jz}(x_z-y_z)(T^{-1})_{jl} 
  }{|T^{-1}(x-y)|^2}\biggr]
\\ \nonumber
&&\   
 -
 \frac{\sum_{l=1}^n\nu_l(x)
 \bigl[
 \nu_l(x) \nu_h(y)  
 -  \nu_h(x)  \nu_l(y) 
 \bigr]}{
 s_n\sqrt{\det a^{(2)} }   |T^{-1}(x-y)|^{n}} 
 \\ \nonumber
&&\quad   
-(2-n)|x-y|^{1-n}A_2\left(\frac{x-y}{|x-y|},|x-y|\right)a^{(2)}\nu(y)
\\ \nonumber
&&\quad   
\times
\sum_{l=1}^n\nu_l(x)\bigl[
\nu_l(x)\frac{x_h-y_h}{|x-y|}-\nu_h(x)\frac{x_l-y_l}{|x-y|}
\bigr]
\\ \nonumber
&&\   
-\sum_{j=1}^n\frac{\partial A_2}{\partial y_j}
\left(\frac{x-y}{|x-y|},|x-y|\right)a^{(2)}\nu(y)|x-y|^{-n}
\\ \nonumber
&&\quad   
\times\sum_{l=1}^n\nu_l(x)
\biggl[
\nu_l(x)\biggl(
\delta_{jh}|x-y|-\frac{(x_j-y_j)(x_h-y_h)}{|x-y|}
\biggr)
\\ \nonumber
&&\   
-\nu_h(x)
\biggl(
\delta_{jl}|x-y|-\frac{(x_j-y_j)(x_l-y_l)}{|x-y|}
\biggr)
\biggr]
\\ \nonumber
&&\   
-\frac{\partial A_2}{\partial r}
\left(\frac{x-y}{|x-y|},|x-y|\right)a^{(2)}\nu(y)
\\ \nonumber
&&\quad   
\times
\sum_{l=1}^n\nu_l(x)
\left[
\nu_l(x)\frac{x_h-y_h}{|x-y|^{n-1}}-\nu_h(x)\frac{x_l-y_l}{|x-y|^{n-1}}
\right]
\\ \nonumber
&&\    
-\sum_{j,s=1}^n
\sum_{l=1}^n\nu_l(x)
\left[
\nu_l(x)
\frac{\partial^2B_1}{\partial x_h\partial x_j}(x-y)
-
\nu_h(x)
\frac{\partial^2B_1}{\partial x_l\partial x_j}(x-y)
\right]
\\ \nonumber
&&\quad   
\times
a_{js}\nu_s(y) \ln |x-y|
\\ \nonumber
&&\   
-DB_1(x-y)a^{(2)}\nu(y) 
\sum_{l=1}^n\nu_l(x)
\left[
\nu_l(x)\frac{x_h-y_h}{|x-y|^2}-\nu_h(x)\frac{x_l-y_l}{|x-y|^2}
\right]
\\ \nonumber
&&\   
-\sum_{j,s=1}^n
\sum_{l=1}^n\nu_l(x)
\left[\nu_l(x)\frac{\partial^2C}{\partial x_h\partial x_j}(x-y)
-
\nu_h(x)\frac{\partial^2C}{\partial x_l\partial x_j}(x-y)
\right]
\\ \nonumber
&&\quad   
\times
a_{js}\nu_s(y) 
\\ \nonumber
&&\   
-\nu(y)^t\cdot a^{(1)} 
\sum_{l=1}^n\nu_l(x)
\left[\nu_l(x)\frac{\partial S_{{\mathbf{a}}} }{\partial x_h}(x-y)
- \nu_h(x)\frac{\partial S_{{\mathbf{a}}} }{\partial x_l}(x-y)
\right] 
 \end{eqnarray}
 for all $(x,y)\in (\partial\Omega)^2\setminus {\mathbb{D}}_{\partial\Omega}$, where we understand that the symbols
 \[
 \frac{\partial A_2}{\partial y_j}\qquad\forall j\in\{1,\dots,n\}
 \]
 denote    partial derivatives of any of the analytic extensions of $A_2$ to an open neighborhood of $(\partial{\mathbb{B}}_n(0,1))\times {\mathbb{R}}$ in ${\mathbb{R}}^{n+1}$.
 \end{lemma}
 For a proof, we refer to \cite{La22d}. Finally, we note that if  $T\in M_{n}({\mathbb{R}})$ is invertible, then
\begin{equation}\label{thm:bmfdl3}
|T^{-1}\xi|\geq |T|^{-1}|\xi|\qquad\forall \xi\in {\mathbb{R}}^n\,,
\end{equation}
where $|T|$ denotes the operator norm of  $T$. (cf.~\textit{e.g.},  \cite[Lem.~3.1 (ii)]{DoLa17} with Dondi).

\section{The maximal function on the boundary of the kernel associated to a positively homogeneous function of degree $-(n-1)$}\label{sec:mafuhomog}
 Let
\begin{equation}\label{volume.klhomog}
{\mathcal{K}}^{0,1}_{-(n-1)} 
\equiv\biggl\{
k\in C^{0,1}_{ {\mathrm{loc}}}({\mathbb{R}}^n\setminus\{0\}):\,
 k\ {\text{is\ positively\ homogeneous\ of \ degree}}\ -(n-1)
\biggr\}\,,
\end{equation}
and
\[
\|k\|_{ {\mathcal{K}}^{0,1}_{-(n-1)}}\equiv \|k\|_{C^{0,1}(\partial{\mathbb{B}}_n(0,1))}
=\sup_{\partial{\mathbb{B}}_n(0,1)}|k|+|k:\,\partial{\mathbb{B}}_n(0,1)|_1
\quad\forall k\in {\mathcal{K}}^{0,1}_{-(n-1)}\,.
\]
Here $C^{0,1}_{ {\mathrm{loc}}}({\mathbb{R}}^n\setminus\{0\})$ denotes the space of functions from ${\mathbb{R}}^n\setminus\{0\}$ to ${\mathbb{C}}$ that are Lipschitz continuous on the compact subsets of ${\mathbb{R}}^n\setminus\{0\}$. One can easily verify that $ \left({\mathcal{K}}^{0,1}_{-(n-1)} , \|\cdot\|_{ {\mathcal{K}}^{0,1}_{-(n-1)}}\right)$ is a Banach space  and we  consider the closed subspace
\begin{equation}\label{volume.k01e0}
{\mathcal{K}}^{0,1}_{-(n-1);o}\equiv \left\{
k\in {\mathcal{K}}^{0,1}_{-(n-1)}:\,k\ {\mathrm{is\ odd}}\right\}
\end{equation}
of ${\mathcal{K}}^{0,1}_{-(n-1)}$. We now to prove that the maximal function on the boundary of a bounded open subset of ${\mathbb{R}}^n$ of class $C^{1,\tilde{\alpha}}$  for some $\tilde{\alpha}\in]0,1]$
and corresponding to a convolution kernel $k(x-y)$, where $k$ is 
an odd locally Lipschitz continuous positively homogeneous  function of degree $-(n-1)$  is bounded. Such a result is well know under stronger assumptions on $k$ 
and one could deduce it from the classical work of Avantaggiati \cite[Thm.~3.1, p.~205]{Av63}, and Miranda \cite[p.~314]{Mi65}, Mitrea, Mitrea and Verdera \cite{MitMitVe16}. Here we provide a direct proof  that improves the regularity of $k$. 
\begin{theorem}\label{thm:mfn-1hbdd}
Let $n\in {\mathbb{N}}$,  $n\geq 2$,   $\tilde{\alpha}\in]0,1]$. Let $\Omega$ be a bounded open subset of ${\mathbb{R}}^n$   of class $C^{1,\tilde{\alpha}}$. Then there exist  $r_{\partial\Omega}$, $\tilde{c}_{\partial\Omega,\tilde{\alpha}}\in]0,+\infty[$ such that 
\begin{equation}\label{thm:mfn-1hbdd1}
\sup_{x\in \partial\Omega}\left|\,
 \int_{(\partial\Omega)\cap 
 [
 {\mathbb{B}}_{n}(x,r) \setminus{\mathbb{B}}_{n}(x,\epsilon)
 ]} k(x-y)  \,d\sigma_{y}
\,
\right|\leq \tilde{c}_{\partial\Omega,\tilde{\alpha}}\left\| k \right\|_{{\mathcal{K}}^{0,1}_{-(n-1);o }}
r^{\tilde{\alpha}}\,,
\end{equation}
for all $\epsilon,r\in]0,r_{\partial\Omega}[$ such that $ \epsilon<r$ and for all $k\in {\mathcal{K}}^{0,1}_{-(n-1);o }$.
\end{theorem}
{\bf Proof.} By the Lemma  of the uniform cylinders, there exist $r_{\partial\Omega}$, $\delta\in]0,1[$ such that for each $\xi\in \partial\Omega$ there exist
$R_\xi\in O_n({\mathbb{R}})$ such that $C(\xi,R_\xi,r_{\partial\Omega},\delta)$ is a coordinate cylinder for $\partial\Omega$ around $\xi$ and that the corresponding function $\gamma_\xi$ satisfies  the conditions
\[
\gamma_\xi(0)=0\,,\quad
D\gamma_\xi(0)=0\qquad \forall\xi\in \partial\Omega\,,\quad
A\equiv\sup_{\xi\in \partial\Omega}\|\gamma_\xi\|_{ C^{1,\tilde{\alpha}}(\overline{{\mathbb{B}}_{n-1}(0,r_{\partial\Omega})}) }<+\infty 
\]
(cf.~\textit{e.g.},   \cite[Lem.~2.63]{DaLaMu21} with Dalla Riva and Musolino). Next we denote by $M_\xi$ the tangent space to $\partial\Omega$ at $\xi$ for all points $\xi$ of $\partial\Omega$. Now let $k$ belong to $ {\mathcal{K}}^{0,1}_{-(n-1);o }$. 
Since $k$ is odd and positively homogeneous  of degree $-(n-1)$, we have
\[
\int_{M_x\cap 
({\mathbb{B}}_{n}(x,r)
\setminus {\mathbb{B}}_{n}(x,\epsilon))
} k(x-y)  \,d\sigma_{y}=0 \qquad\forall x\in \partial\Omega
\]
for all $\epsilon$, $r\in ]0,+\infty[$ such that $\epsilon<r$. Then the triangular inequality implies that 
\begin{eqnarray}\label{thm:mfn-1hbdd2}
\lefteqn{
\biggl|\biggr.
\int_{(\partial\Omega)\cap 
[{\mathbb{B}}_{n}(x,r)
\setminus {\mathbb{B}}_{n}(x,\epsilon)]
} k(x-y)  \,d\sigma_{y}\biggl.\biggr|
}
\\ \nonumber
&& 
\leq
\biggl|\biggr.
\int_{(\partial\Omega)\cap 
[{\mathbb{B}}_{n}(x,r)
\setminus {\mathbb{B}}_{n}(x,\epsilon)]
} k(x-y)  \,d\sigma_{y}
\\ \nonumber
&&  \quad
-
\int_{M_x\cap 
[{\mathbb{B}}_{n}(x,r)
\setminus {\mathbb{B}}_{n}(x,\epsilon)]
} k(x-y)  \,d\sigma_{y}\biggl.\biggr|
\\ \nonumber
&&  \quad
+\biggl|\biggr.
\int_{M_x\cap 
[{\mathbb{B}}_{n}(x,r)
\setminus {\mathbb{B}}_{n}(x,\epsilon)]
} k(x-y)  \,d\sigma_{y}\biggl.\biggr|
\\ \nonumber
&& 
=\biggl|\biggr.
\int_{\{
\eta\in {\mathbb{B}}_{n-1}(0,r_{\partial\Omega}):\,\epsilon^2\leq|\eta|^2+|\gamma_x(\eta)|^2<r^2
\}}
k\circ R_x^t ((-\eta,-\gamma_x(\eta))^t)
\\ \nonumber
&&  \quad
\times
\sqrt{1+|D\gamma_x(\eta)|^2}\,d\eta
\\ \nonumber
&&\quad
-
\int_{\{
\eta\in {\mathbb{B}}_{n-1}(0,r_{\partial\Omega}):\,\epsilon^2\leq|\eta|^2 <r^2
\}}
k\circ R_x^t ((-\eta,0)^t)
 \,d\eta
\biggl.\biggr|\leq I_1+I_2
\end{eqnarray}
for all $x\in  \partial\Omega$ and  $\epsilon$, $r\in ]0,r_{\partial\Omega}[$ such that $\epsilon<r$, where $I_1$ and $I_2$ are as follows
\begin{eqnarray*}
\lefteqn{I_1
\equiv
  \biggl|\biggr.
  \int_{\{
\eta\in {\mathbb{B}}_{n-1}(0,r_{\partial\Omega}):\,\epsilon^2\leq|\eta|^2+|\gamma_x(\eta)|^2<r^2
\}}
k\circ R_x^t ((-\eta,-\gamma_x(\eta))^t)
\sqrt{1+|D\gamma_x(\eta)|^2}\,d\eta
}
\\ \nonumber
&&\qquad 
-
 \int_{\{
\eta\in {\mathbb{B}}_{n-1}(0,r_{\partial\Omega}):\,\epsilon^2\leq|\eta|^2 <r^2
\}}
k\circ R_x^t ((-\eta,-\gamma_x(\eta))^t)
\sqrt{1+|D\gamma_x(\eta)|^2}\,d\eta
  \biggl.\biggr|\,,
\\ \nonumber
\lefteqn{I_2\equiv
\biggl|\biggr.
 \int_{\{
\eta\in {\mathbb{B}}_{n-1}(0,r_{\partial\Omega}):\,\epsilon^2\leq|\eta|^2 <r^2
\}}
k\circ R_x^t ((-\eta,-\gamma_x(\eta))^t)
\sqrt{1+|D\gamma_x(\eta)|^2}\,d\eta
}
\\ \nonumber
&&  \qquad
-
\int_{\{
\eta\in {\mathbb{B}}_{n-1}(0,r_{\partial\Omega}):\,\epsilon^2\leq|\eta|^2 <r^2
\}}
k\circ R_x^t ((-\eta,0)^t)
\,d\eta\biggl.\biggr|\,. 
\end{eqnarray*}
In order to estimate $I_1$ and $I_2$ we note that  the Taylor formula   for the function $\gamma_x$ and   equalities $\gamma_x(0)=0$ and $D\gamma_x(0)=0$ imply that
\begin{equation}\label{thm:mfn-1hbdd6}
|\gamma_x(\eta)|\leq
 | D \gamma_x |_{\tilde{\alpha}}
   |\eta|^{\tilde{\alpha}+1}
   \leq A |\eta|^{\tilde{\alpha}+1}
    \qquad\forall \eta\in \overline{{\mathbb{B}}_{n-1}(0,r_{\partial\Omega})} 
\end{equation}
(cf.~\textit{e.g.}, \cite[Lem.~2.58]{DaLaMu21}). Then we set
\[
a_x(r)\equiv
 \sup	 
\left\{
\frac{|\gamma_x(\eta)|}{|\eta|}:\,\eta\in \overline{{\mathbb{B}}_{n-1}(0,r)}\setminus\{0\}
\right\}\qquad\forall r\in]0,r_{\partial\Omega}[ 
\]
and we note that
\[
a_x(r)\leq | D \gamma_x |_{\tilde{\alpha}} r^{\tilde{\alpha}} 
\leq A  r^{\tilde{\alpha}} 
\qquad\forall r\in]0,r_{\partial\Omega}[ 
\]
and that
\begin{equation}\label{thm:mfn-1hbdd7}
|\eta|^2+|\gamma_x(\eta)|^2\leq
|\eta|^2 (1+a_x^2(r) ) \qquad\forall \eta\in \overline{{\mathbb{B}}_{n-1}(0,r)}
\end{equation}
for all $r\in]0,r_{\partial\Omega}[$ and $x\in  \partial\Omega $. In order to estimate $I_1$, we plan to apply Lemma \ref{lem:pvcom} of the Appendix. By inequality (\ref{thm:mfn-1hbdd6}), we have
\begin{eqnarray*}
\lefteqn{
c_k\equiv \sup_{\eta\in {\mathbb{B}}_{n-1}(0,r_{\partial\Omega})\setminus\{0\}} 
|\eta|^{(n-1)}\biggl|
k\circ R_x^t ((-\eta,-\gamma_x(\eta))^t)
\sqrt{1+|D\gamma_x(\eta)|^2} \biggr|
}
\\ \nonumber
&&\qquad\qquad\qquad\qquad\qquad\qquad
\leq
 \left(\sup_{\partial {\mathbb{B}}_{n-1}(0,1)}|k|\right)\sqrt{1+A^2} 
 \,.
\end{eqnarray*}
Then Lemma \ref{lem:pvcom} implies that
\begin{eqnarray} \nonumber
&&
I_1=\biggl|\biggr.
  \int_{\{
\eta\in {\mathbb{B}}_{n-1}(0,r):\,\epsilon^2\leq|\eta|^2+|\gamma_x(\eta)|^2 
\}}
k\circ R_x^t ((-\eta,-\gamma_x(\eta))^t)
\sqrt{1+|D\gamma_x(\eta)|^2}\,d\eta
\\ \label{thm:mfn-1hbdd7a}
&&\quad 
-
 \int_{\{
\eta\in {\mathbb{B}}_{n-1}(0,r):\,\epsilon^2\leq|\eta|^2   
\}}
k\circ R_x^t ((-\eta,-\gamma_x(\eta))^t)
\sqrt{1+|D\gamma_x(\eta)|^2}\,d\eta
  \biggl.\biggr|
\\ \nonumber
&&\quad\quad  
  \leq \left(\sup_{\partial {\mathbb{B}}_{n-1}(0,1)}|k|\right)\sqrt{1+A^2}s_{n-1}\log(1+a_x(\epsilon))
 \\ \nonumber
&&\quad\quad  
  \leq \left(\sup_{\partial {\mathbb{B}}_{n-1}(0,1)}|k|\right)\sqrt{1+A^2}s_{n-1}a_x(\epsilon)
 \leq \left(\sup_{\partial {\mathbb{B}}_{n-1}(0,1)}|k|\right)\sqrt{1+A^2}s_{n-1}
 A  r^{\tilde{\alpha}}  
\end{eqnarray}
for all $x\in  \partial\Omega$ and  $\epsilon$, $r\in ]0,r_{\partial\Omega}[$ such that $\epsilon<r$.  (cf.~inequality (\ref{lem:pvcom1b}) with $a=0$).  We now turn  to estimate $I_2$. Since $k$ is positively homogeneous of degree $-(n-1)$, the inequality of Cialdea \cite[VIII, p.~47]{Ci95} (see also 
\cite[Lem.~4.14]{DaLaMu21}, case $\alpha=1$) implies that
\begin{eqnarray}\label{thm:mfn-1hbdd8}
\lefteqn{
|k\circ R_x^t ((-\eta,-\gamma_x(\eta))^t)-k\circ R_x^t ((-\eta,0)^t)|
}
\\ \nonumber
&&\qquad
\leq
2n \|k\|_{C^{0,1}(\partial {\mathbb{B}}_{n}(0,1))}|R_x^t ((0,-\gamma_x(\eta))^t)|
\\ \nonumber
&&  \qquad\quad
\times\min\{
|R_x^t ((-\eta,-\gamma_x(\eta))^t)|,
|R_x^t ((-\eta,0)^t)|
\}^{-n}
\\ \nonumber
&&  \qquad
\leq 2n \|k\|_{C^{0,1}(\partial {\mathbb{B}}_{n}(0,1))}
|\gamma_x(\eta)|\,|\eta|^{-n}
\\ \nonumber
&&  \qquad
\leq 2n \|k\|_{C^{0,1}(\partial {\mathbb{B}}_{n}(0,1))}
 | D \gamma_x |_{\tilde{\alpha}}\,|\eta|^{-(n-1)+\tilde{\alpha}}
\ \ \forall\eta\in  {\mathbb{B}}_{n-1}(0,r_{\partial\Omega})\setminus\{0\}\,.
\end{eqnarray}
Since the function $(1+t^2)^{1/2}$ is Lipschitz continuous in the variable $t\in {\mathbb{R}}$ with Lipschitz constant equal to $1$, we have
\begin{equation}\label{thm:mfn-1hbdd9}
\left|
\sqrt{1+|D\gamma_x(\eta)|^2}-1
\right|\leq |\,|D\gamma_x(\eta)|-|D\gamma_x(0)|\,|\leq  | D \gamma_x |_{\tilde{\alpha}}|\eta|^{\tilde{\alpha}}  
\end{equation}
for all $\eta\in  {\mathbb{B}}_{n-1}(0,r_{\partial\Omega})$. Then inequalities (\ref{thm:mfn-1hbdd8}) and (\ref{thm:mfn-1hbdd9}) imply that
\begin{eqnarray}\label{thm:mfn-1hbdd10}
\lefteqn{
I_2
\leq
\biggl|\biggr.
\int_{\{
\eta\in {\mathbb{B}}_{n-1}(0,r ):\,\epsilon^2<|\eta|^2 
\}}
[k\circ R_x^t ((-\eta,-\gamma_x(\eta))^t)
}
\\ \nonumber
&&  \quad 
-
k\circ R_x^t ((-\eta,0)^t)]
\sqrt{1+|D\gamma_x(\eta)|^2}\,d\eta
\biggl.\biggr|
\\ \nonumber
&&\quad
+
\int_{\{
\eta\in {\mathbb{B}}_{n-1}(0,r):\,\epsilon^2<|\eta|^2 
\}}
\left| k\circ R_x^t ((-\eta,0)^t) 
[\sqrt{1+|D\gamma_x(\eta)|^2}-1]
\right|\,d\eta  
\\ \nonumber
&& 
\leq
2n \|k\|_{C^{0,1}(\partial {\mathbb{B}}_{n}(0,1))}
\sqrt{1+A^2} | D \gamma_x |_{\tilde{\alpha}}
\\ \nonumber
&&  \quad
\times\int_{\{
\eta\in {\mathbb{B}}_{n-1}(0,r ):\,\epsilon^2<|\eta|^2 
\}}|\eta|^{-(n-1)+\tilde{\alpha}}\,d\eta
\\ \nonumber
&&\quad
+\sup_{\partial {\mathbb{B}}_{n-1}(0,1)}|k|
 | D \gamma_x |_{\tilde{\alpha}}
 \int_{\{
\eta\in {\mathbb{B}}_{n-1}(0,r):\,\epsilon^2<|\eta|^2 
\}}|\eta|^{-(n-1)+\tilde{\alpha}}\,d\eta
\\ \nonumber
&& 
\leq
(2n\sqrt{1+A^2}+1) \|k\|_{C^{0,1}(\partial {\mathbb{B}}_{n}(0,1))}\|\gamma_x\|_{
C^{1,\tilde{\alpha}}(\overline{ {\mathbb{B}}_{n-1}(0,r_{\partial\Omega})})
}
\\ \nonumber
&&  \quad
\times\frac{s_{n-1}}{\tilde{\alpha}} ( r^{\tilde{\alpha}}-\epsilon^{\tilde{\alpha}})
\\ \nonumber
&& 
\leq
A(2n\sqrt{1+A^2}+1) 
\frac{s_{n-1}}{\tilde{\alpha}}   r^{\tilde{\alpha}}
\|k\|_{C^{0,1}(\partial {\mathbb{B}}_{n}(0,1))}  
\qquad\forall r\in]0,r_{\partial\Omega}[\,.
\end{eqnarray}
By inequalities    (\ref{thm:mfn-1hbdd2}),  (\ref{thm:mfn-1hbdd7a}),  (\ref{thm:mfn-1hbdd10}), we conclude that the constant $\tilde{c}_{\partial\Omega,\tilde{\alpha}}$ of the inequality of the statement exists and the proof is complete.\hfill  $\Box$ 

\vspace{\baselineskip}

Then we have the following immediate consequence of Theorem \ref{thm:mfn-1hbdd}.
\begin{theorem}\label{lem:maolg}
 Let $n\in {\mathbb{N}}$,  $n\geq 2$,   $\tilde{\alpha}\in]0,1]$.  Let $\Omega$ be a bounded open subset of ${\mathbb{R}}^n$ of class $C^{1,\tilde{\alpha}}$. Then there exists   $c^*_{\partial\Omega,\tilde{\alpha}}\in]0,+\infty[$ such that 
\begin{equation}\label{lem:maolg1}
\sup_{x\in \partial\Omega}
\sup_{\epsilon\in ]0,+\infty[}
\left|
\int_{(\partial\Omega)\setminus{\mathbb{B}}_{n}(x,\epsilon)}
k(x-y)\,d\sigma_y
\right|
\leq c^*_{\partial\Omega,\tilde{\alpha}}\left\| k \right\|_{{\mathcal{K}}^{0,1}_{-(n-1);o }}
\ \ \forall k\in {\mathcal{K}}^{0,1}_{-(n-1);o }\,.
\end{equation}
\end{theorem}
Indeed, $|k(x-y)|\leq \sup_{\partial {\mathbb{B}}_{n-1}(0,1)}|k|(r_{\partial\Omega}/2)^{-(n-1)}$ if $(x,y)\in (\partial\Omega)^2$ and $ |x-y|\geq r_{\partial\Omega}/2$, where $r_{\partial\Omega}$ is as in Theorem \ref{thm:mfn-1hbdd}.

  \subsection{Boundedness of the maximal function of   the tangential gradient of the kernel of the double layer potential}
  In order to prove inequality (\ref{thm:dllreggenn1}), we first estimate the maximal function of the kernel of the first order derivatives of the fundamental solution $S_{{\mathbf{a}}}$. 
\begin{theorem}\label{thm:maggnfs}
  Let $n\in {\mathbb{N}}\setminus\{0,1\}$. 
 Let ${\mathbf{a}}$ be as in (\ref{introd0}), (\ref{ellip}), (\ref{symr}).  Let $S_{ {\mathbf{a}} }$ be a fundamental solution of $P[{\mathbf{a}},D]$.   Let $\alpha\in]0,1[$.  Let $\Omega$ be a bounded open subset of ${\mathbb{R}}^n$  of class $C^{1, \alpha }$. Then \begin{equation}\label{thm:maggnfs1}
\sup_{x\in \partial\Omega}
\sup_{r\in ]0,+\infty[}
\left|\,
\int_{(\partial\Omega)\setminus{\mathbb{B}}_{n}(x,r)}
\frac{\partial S_{{\mathbf{a}}}}{\partial x_h}(x-y)\,d\sigma_y
\right|
<+\infty
\end{equation}
for all $h\in\{1,\dots,n\}$.
\end{theorem}
{\bf Proof.}  By formula (\ref{grafun1}) and by the known inequalities 
\begin{equation}\label{thm:maggnfs2}
\sup_{x\in \partial\Omega}\int_{\partial\Omega}|x-y|^{-\gamma}\,d\sigma_y<+\infty\,,
\quad
\sup_{x\in \partial\Omega}\int_{\partial\Omega}
|\ln |x-y||
\,d\sigma_y<+\infty
\end{equation}
for $\gamma\in]-\infty,(n-1)[$ (cf.~\textit{e.g.},  \cite[Lem.~3.5]{DoLa17}),  we have
\begin{eqnarray*}
\lefteqn{
\sup_{x\in \partial\Omega}
\sup_{r\in ]0,+\infty[}
\biggl| \,
\int_{(\partial\Omega)\setminus{\mathbb{B}}_{n}(x,r)}
DS_{ {\mathbf{a}} }(x-y)
}
\\ \nonumber
&&  \quad
-\frac{1}{ s_{n}\sqrt{\det a^{(2)} } }
|T^{-1}(x-y)|^{-n}(x-y)^{t}(a^{(2)})^{-1}\,d\sigma_y
\biggr|
\\ \nonumber
&&\qquad
\leq 
\sup_{x\in \partial\Omega}
\int_{\partial\Omega} \biggl|	 
 |	 
x-y|^{2-n}A_{2}(\frac{x-y}{|x-y|},|x-y|)
\\ \nonumber
&&  \quad
+DB_{1}(x-y)\ln |x-y|+DC(x-y)
\biggr|\,d\sigma_y<+\infty\,.
\end{eqnarray*}
Since the function $\frac{1}{ s_{n}\sqrt{\det a^{(2)} } }
|T^{-1}\xi|^{-n}\xi^{t}(a^{(2)})^{-1}$ is positively homogeneous of degree $-(n-1)$, Theorem \ref{lem:maolg} implies  that 
\[
\sup_{x\in \partial\Omega}
\sup_{r\in ]0,+\infty[}
\left|\,
\int_{(\partial\Omega)\setminus{\mathbb{B}}_{n}(x,r)}
\frac{1}{ s_{n}\sqrt{\det a^{(2)} } }
|T^{-1}(x-y)|^{-n}(x-y)^{t}(a^{(2)})^{-1}\,d\sigma_y
\right|
\]
is finite. Then the above  inequality implies the validity of the statement.
\hfill  $\Box$ 

\vspace{\baselineskip}

 Next we prove the following intermediate statement.
\begin{proposition}\label{prop:mftgdlgech}
 Let $n\in {\mathbb{N}}\setminus\{0,1\}$,  $\beta\in]0,1[$.  Let $\Omega$ be a bounded open subset of ${\mathbb{R}}^{n}$ of class $C^{1,1}$.  
Let ${\mathbf{a}}$ be as in (\ref{introd0}), (\ref{ellip}), (\ref{symr}).  Let $T\in M_{n}({\mathbb{R}})$  be as in (\ref{prop:ourfs0}).  Let $S_{ {\mathbf{a}} }$ be a fundamental solution of $P[{\mathbf{a}},D]$. Then condition (\ref{thm:dllreggenn1}) is equivalent to the following condition
\begin{eqnarray}\label{thm:mftgdlgx1}
&&
\sup_{x\in \partial\Omega}\sup_{r\in ]0,+\infty[}
\biggl\vert
\int_{(\partial\Omega)\setminus  {\mathbb{B}}_n(x,r)}
\\ \nonumber
&&\qquad
{\mathrm{grad}}_{\partial\Omega,x}\left(
-\frac{1}{s_{n}\sqrt{\det a^{(2)} }}|T^{-1}(x-y)|^{-n} (x-y)^{t}\nu (y)
\right)\,d\sigma_y
\biggr\vert <+\infty\,.
\end{eqnarray}
\end{proposition}
 {\bf Proof.} By formulas (\ref{eq:boest3}), (\ref{lem:tgdlgen1}), we have
\begin{eqnarray*} 
\lefteqn{
  {\mathrm{grad}}_{\partial\Omega,x}\overline{B^{*}_{\Omega,y}}\left(S_{{\mathbf{a}}}(x-y)\right) 
}
\\ \nonumber
&&\qquad\qquad
= {\mathrm{grad}}_{\partial\Omega,x}\left(
-\frac{1}{s_{n}\sqrt{\det a^{(2)} }}|T^{-1}(x-y)|^{-n} (x-y)^{t}\nu (y)
\right) 
\\ \nonumber
&&\qquad\qquad\quad
+\sum_{j=3}^9\left(J_{jh}(x,r)\right)_{h=1,\dots,n} \,,
\end{eqnarray*}
for all $ x\in\partial\Omega, r\in]0,+\infty[$, 
where $J_{jh}(x,r)$ denotes  the  $j$-th addendum of the right hand side of formula (\ref{lem:tgdlgen1}) for the $h$-th component of 
\[
{\mathrm{grad}}_{\partial\Omega,x}\overline{B^{*}_{\Omega,y}}\left(S_{{\mathbf{a}}}(x-y)\right)\,.
\]
 Then the equivalence of (\ref{thm:dllreggenn1}) and (\ref{thm:mftgdlgx1}) follows by the following inequalities
\begin{eqnarray}\label{thm:mftgdlgech2c}
\lefteqn{
\sup_{x\in\partial\Omega}\sup_{r\in]0,+\infty[} 
\biggl|\int_{(\partial\Omega)\setminus  {\mathbb{B}}_n(x,r)}J_{jh}(x,r)\,d\sigma_y\biggr|<+\infty
}
\\ \nonumber
&&\qquad\qquad\qquad\qquad\quad\forall j\in\{3,\dots,9\}\,,
h\in\{1,\dots,n\}\,.
\end{eqnarray}
We now fix $h\in\{1,\dots,n\}$ and we turn to estimate 
 the integral of	 
$J_{jh}(x,r)$
 on	$(\partial\Omega)\setminus  {\mathbb{B}}_n(x,r)$	 
for each $j\in\{3,\dots,9\}$. We first look at  the   third addendum $J_{3h}(x,r)$ in the right hand side of formula (\ref{lem:tgdlgen1}). Since $A_2$ is analytic, we have
\begin{eqnarray*}
\lefteqn{
 A_2\left(\frac{x-y}{|x-y|},|x-y|\right)
 =\int_0^1\frac{\partial A_2}{\partial r}\left(\frac{x-y}{|x-y|},s|x-y|\right)\,ds|x-y|  
 }
\\ \nonumber
&&\qquad\qquad\qquad\qquad\qquad\qquad\qquad\qquad
+A_2\left(\frac{x-y}{|x-y|},0\right)
\end{eqnarray*} 
for all $(x,y)\in  (\partial\Omega)^2 \setminus {\mathbb{D}}_{ \partial\Omega }$. 
Then we note that 
\[
a_{2}'\equiv\sup_{(x,y)\in (\partial\Omega)^2\setminus {\mathbb{D}}_{\partial\Omega}}
\left| \int_0^1\frac{\partial A_2}{\partial r}\left(\frac{x-y}{|x-y|},s|x-y|\right)\,ds \right|<+\infty
\]
and that the function $A_2\left(\frac{\xi}{|\xi|},0\right)$ is positively homogeneous of degree $0$ in the variable $\xi\in {\mathbb{R}}^n$. By  Proposition \ref{prop:grafun}   $A_2(\cdot,0)$ is even. Then the function 
\[
 |\xi|^{1-n}A_2\left(\frac{\xi}{|\xi|},0\right)\frac{\xi_l}{|\xi|}\qquad\forall \xi\in 
{\mathbb{R}}^n\setminus\{0\}
\]
is positively homogeneous of degree $-(n-1)$  and odd for all $l\in\{1,\dots,n\}$ and Theorem \ref{lem:maolg} on the estimate of the maximal function associated to an odd positively homogeneous kernel of degree $-(n-1)$ implies that
\begin{eqnarray}\label{thm:mftgdlgech9}
 \lefteqn{
 \biggl|\int_{\partial\Omega \setminus{\mathbb{B}}_n(x,r)  }
 |x-y|^{1-n}A_2\left(\frac{x-y}{|x-y|},|x-y|\right)a^{(2)}\nu(y)
 }
\\ \nonumber
&&\qquad
 \times\sum_{l=1}^n\nu_l(x)\bigl[
\nu_l(x)\frac{x_h-y_h}{|x-y|}-\nu_h(x)\frac{x_l-y_l}{|x-y|}
\bigr]\,d\sigma_y\biggr|
\\ \nonumber
&&\qquad
\leq \int_{\partial\Omega}|x-y|^{2-n}\left|
\int_0^1\frac{\partial A_2}{\partial r}\left(\frac{x-y}{|x-y|},s|x-y|\right)\,ds\right||a^{(2)}|2n\,d\sigma_y
\\ \nonumber
&&\qquad\quad
+\biggl|\int_{\partial\Omega}
 |x-y|^{1-n}A_2\left(\frac{x-y}{|x-y|},0\right)a^{(2)}(\nu(y)-\nu(x))
\\ \nonumber
&&\qquad
 \times
 \sum_{l=1}^n\nu_l(x)\bigl[
\nu_l(x)\frac{x_h-y_h}{|x-y|}-\nu_h(x)\frac{x_l-y_l}{|x-y|}
\bigr]\,d\sigma_y\biggr|
\\ \nonumber
&&\qquad\quad
+\biggl|\int_{\partial\Omega \setminus{\mathbb{B}}_n(x,r)  }
 |x-y|^{1-n}A_2\left(\frac{x-y}{|x-y|},0\right)a^{(2)}\nu(x)
\\ \nonumber
&&\qquad
\times \sum_{l=1}^n\nu_l(x)\bigl[
\nu_l(x)\frac{x_h-y_h}{|x-y|}-\nu_h(x)\frac{x_l-y_l}{|x-y|}
\bigr]\,d\sigma_y\biggr|
\\ \nonumber
&& \qquad
\leq 
a_{2}'\sup_{\xi\in\partial\Omega}  \int_{\partial\Omega}|\xi-y|^{2-n}\,d\sigma_y|a^{(2)}|2n
\\ \nonumber
&&\qquad\quad
+\sup_{\eta\in\partial {\mathbb{B}}_n(0,1)}|A_2(\eta,0)|\,| a^{(2)}|2n{\mathrm{Lip}} (\nu) \sup_{\xi\in\partial\Omega}  \int_{\partial\Omega}|\xi-y|^{2-n}\,d\sigma_y
\\ \nonumber
&&\qquad \quad
+
2n| a^{(2)}| \sup_{	l\in\{1,\dots,n\}	}\sup_{\xi\in\partial\Omega}\sup_{\rho\in]0,+\infty[}
\biggl|
\int_{\partial\Omega \setminus{\mathbb{B}}_n(\xi,\rho)  }|\xi-y|^{1-n} 
\\ \nonumber
&&\qquad
 \times A_2\left(\frac{\xi-y}{|\xi-y|},0\right) 
\frac{\xi_l-y_l}{|\xi-y|}
\,d\sigma_y
\biggr|<+\infty\,,
\end{eqnarray}
for all $x\in \partial\Omega$ and $r\in]0,+\infty[$ (see also (\ref{thm:maggnfs2})). We now look at the fourth addendum $J_{4h}(x,r)$ in the right hand side of formula (\ref{lem:tgdlgen1}). Here we note that   formula (\ref{lem:tgdlgen1}) holds when the symbol 
 $\frac{\partial A_2}{\partial y_j}$
denotes   partial derivative with respect to $y_j$ of any of the analytic extensions of $A_2$ to an open neighborhood of the domain of $A_2$. Here we choose the analytic extension of $A_2$ defined by 
\[
\tilde{A}_2(x,r)\equiv A_2(\frac{x}{|x|},r)\qquad\forall (x,r)\in ({\mathbb{R}}^n\setminus\{0\}) \times{\mathbb{R}}
\]
(see  Lemma \ref{lem:anexsph})).  By Proposition \ref{prop:grafun}, $A_2(\cdot,0)$ is even. Then  
$\tilde{A}_2(\cdot,0)$ is also even and thus 
$\frac{\partial \tilde{A}_2}{\partial y_j}(\cdot,0)$ is odd for all $j\in \{0,\dots,n\}$. Moreover, we have
\begin{eqnarray*}
 \lefteqn{
 \frac{\partial \tilde{A}_2}{\partial y_j}\left(\frac{x-y}{|x-y|},|x-y|\right)
 =\int_0^1\frac{\partial^2 \tilde{A}_2}{\partial y_j\partial r}\left(\frac{x-y}{|x-y|},s|x-y|\right)\,ds|x-y|
 }
\\ \nonumber
&&\qquad\qquad\qquad\qquad 
 +\frac{\partial \tilde{A}_2}{\partial y_j}\left(\frac{x-y}{|x-y|},0\right)\qquad\forall (x,y)\in ({\mathbb{R}}^n\times {\mathbb{R}}^n)\setminus {\mathbb{D}}_{{\mathbb{R}}^n}\,,
\end{eqnarray*}
for all $j\in\{1,\dots,n\}$. Then we note that 
\[
a_{2}''\equiv \sup_{ j\in\{1,\dots,n\} }
\sup_{(x,y)\in (\partial\Omega)^2\setminus {\mathbb{D}}_{\partial\Omega}}
\left| \int_0^1\frac{\partial^2 \tilde{A}_2}{\partial y_j\partial r}\left(\frac{x-y}{|x-y|},s|x-y|\right)\,ds \right|<+\infty
\]
and that the function $\frac{\partial \tilde{A}_2}{\partial y_j}\left(\frac{\xi}{|\xi|},0\right)$ is odd and positively homogeneous of degree $0$  in the variable $\xi\in {\mathbb{R}}^n$ for each $j\in\{1,\dots,n\}$. Then the function 
\[
\frac{\partial \tilde{A}_2}{\partial y_j}
\left(\frac{\xi}{|\xi|},0\right) |\xi|^{-n}\biggl(
\delta_{jl}|\xi|-\frac{\xi_j\xi_l}{|\xi|}
\biggr)
\qquad\forall \xi\in 
{\mathbb{R}}^n\setminus\{0\}
\]
is positively homogeneous of degree $-(n-1)$  and odd  for all $j,l\in\{1,\dots,n\}$ and Theorem \ref{lem:maolg} on the estimate of the maximal function associated to an odd positively homogeneous kernel of degree $-(n-1)$ implies that
\begin{eqnarray}\label{thm:mftgdlgech10}
  && \biggl|\int_{\partial\Omega \setminus{\mathbb{B}}_n(x,r)  }
\sum_{j=1}^n\frac{\partial \tilde{A}_2}{\partial y_j}
\left(\frac{x-y}{|x-y|},|x-y|\right) 
\\ \nonumber
&& 
\times a^{(2)}\nu(y)|x-y|^{-n}
\sum_{l=1}^n\nu_l(x)
\biggl[
\nu_l(x)\biggl(
\delta_{jh}|x-y|-\frac{(x_j-y_j)(x_h-y_h)}{|x-y|}
\biggr)
\\ \nonumber
&&\quad 
-\nu_h(x)
\biggl(
\delta_{jl}|x-y|-\frac{(x_j-y_j)(x_l-y_l)}{|x-y|}
\biggr)
\biggr]\,d\sigma_y\biggr|
\\ \nonumber
&&\leq\sum_{j=1}^n\int_{\partial\Omega}\left| \int_0^1\frac{\partial^2 \tilde{A}_2}{\partial y_j\partial r}\left(\frac{x-y}{|x-y|},s|x-y|\right)\,ds\right|
\\ \nonumber
&&\qquad
\times  |a^{(2)}|\,|\nu(y)|\,|x-y|^{2-n}4n \,d\sigma_y
\\ \nonumber
&&\quad+
  \int_{\partial\Omega }
\sum_{j=1}^n\biggl|\frac{\partial \tilde{A}_2}{\partial y_j}
\left(\frac{x-y}{|x-y|},0\right)\biggr|\,|a^{(2)}|\,|\nu(y)-\nu(y)|
\\ \nonumber
&&\quad 
\times |x-y|^{-n}4n|x-y|\,d\sigma_y
\\ \nonumber
&&\quad 
+
 \biggl|\int_{\partial\Omega \setminus{\mathbb{B}}_n(x,r)  }
\sum_{j=1}^n\frac{\partial \tilde{A}_2}{\partial y_j}
\left(\frac{x-y}{|x-y|},0\right)a^{(2)}\nu(x)|x-y|^{-n}
\\ \nonumber
&&\quad
\times
\sum_{l=1}^n\nu_l(x)
\biggl[
\nu_l(x)\biggl(
\delta_{jh}|x-y|-\frac{(x_j-y_j)(x_h-y_h)}{|x-y|}
\biggr)
\\ \nonumber
&&\quad 
-\nu_h(x)
\biggl(
\delta_{jl}|x-y|-\frac{(x_j-y_j)(x_l-y_l)}{|x-y|}
\biggr)
\biggr]\,d\sigma_y\biggr|
\\ \nonumber
&& 
\leq 4n^2a_{2}''|a^{(2)}|\sup_{\xi\in\partial\Omega}\int_{\partial\Omega}|\xi-y|^{2-n}\,d\sigma_y
\\ \nonumber
&&\qquad
+\sum_{j=1}^n\sup_{\eta\in\partial {\mathbb{B}}_n(0,1)}
\biggl|\frac{\partial \tilde{A}_2}{\partial y_j}
\left(\eta,0\right)\biggr|\,|a^{(2)}|\,{\mathrm{Lip}} (\nu)  \sup_{\xi\in\partial\Omega}\int_{\partial\Omega}|\xi-y|^{2-n}\,d\sigma_y4n
\\ \nonumber
&&\quad
+2n^2|a^{(2)}|\sup_{j,l\in\{1,\dots,n\}}\sup_{\xi\in\partial\Omega}\sup_{\rho\in]0,+\infty[}
\biggl|\int_{\partial\Omega \setminus{\mathbb{B}}_n(\xi,\rho)  }
\frac{\partial \tilde{A}_2}{\partial y_j}
\left(\frac{\xi-y}{|\xi-y|},0\right) 
\\ \nonumber
&&\qquad
\times|\xi-y|^{-n}
\biggl(
\delta_{jl}|\xi-y|-\frac{(\xi_j-y_j)(\xi_l-y_l)}{|\xi-y|}
\biggr)\,d\sigma_y
\biggr|<+\infty
\end{eqnarray}
for all $x\in \partial\Omega$ and $r\in]0,+\infty[$ (see also (\ref{thm:maggnfs2})). We now look at the fifth addendum $J_{5h}(x,r)$ in the right hand side of formula (\ref{lem:tgdlgen1}). Since $\frac{\partial A_2}{\partial r}
\left(\frac{x-y}{|x-y|},|x-y|\right)$ is bounded in $(x,y)\in(\partial\Omega)^2\setminus {\mathbb{D}}_{\partial\Omega}$, we have
\begin{eqnarray}\label{thm:mftgdlgech11}
 &&\quad
 \biggl|\int_{\partial\Omega \setminus{\mathbb{B}}_n(x,r)  }
-\frac{\partial A_2}{\partial r}
\left(\frac{x-y}{|x-y|},|x-y|\right)a^{(2)}\nu(y)
\\ \nonumber
&&\quad 
\times\sum_{l=1}^n\nu_l(x)
\left[
\nu_l(x)\frac{x_h-y_h}{|x-y|^{n-1}}-\nu_h(x)\frac{x_l-y_l}{|x-y|^{n-1}}
\right]\,d\sigma_y\biggr|
\\ \nonumber
&&
\leq\sup_{(x,y)\in(\partial\Omega)^2\setminus {\mathbb{D}}_{\partial\Omega}}
\biggl|
\frac{\partial A_2}{\partial r}
\left(\frac{x-y}{|x-y|},|x-y|\right)\biggr| |a^{(2)}|2n
\sup_{\xi\in\partial\Omega}\int_{\partial\Omega}|\xi-y|^{2-n}\,d\sigma_y
<+\infty
\end{eqnarray}
for all $x\in \partial\Omega$ and $r\in]0,+\infty[$  (see also (\ref{thm:maggnfs2})). We now look at the sixth addendum $J_{6h}(x,r)$ in the right hand side of formula (\ref{lem:tgdlgen1}). Since $B_1$ is analytic in ${\mathbb{R}}^n$, we have
\begin{eqnarray}\nonumber
&&
 \biggl|\int_{\partial\Omega \setminus{\mathbb{B}}_n(x,r)  }
\sum_{j,s=1}^n
\sum_{l=1}^n\nu_l(x)
\left[
\nu_l(x)
\frac{\partial^2B_1}{\partial x_h\partial x_j}(x-y)
-
\nu_h(x)
\frac{\partial^2B_1}{\partial x_l\partial x_j}(x-y)
\right]
\\  \label{thm:mftgdlgech12}
&&\quad 
\times a_{js}\nu_s(y) \ln |x-y|\,d\sigma_y\biggr|
\\ \nonumber
&& 
\leq 2 n^3|a^{(2)}| \sup_{j,l\in\{1,\dots,n\}}\sup_{\xi\in(\partial\Omega-\partial\Omega)} 
\left|\frac{\partial^2B_1}{\partial x_l\partial x_j}(\xi)\right| 
\sup_{\xi\in \partial\Omega}\int_{\partial\Omega}|\ln |\xi-y||\,d\sigma_y<+\infty
\end{eqnarray}
for all $x\in \partial\Omega$ and $r\in]0,+\infty[$ (see also (\ref{thm:maggnfs2})). We now look at the seventh addendum $J_{7h}(x,r)$ in the right hand side of formula (\ref{lem:tgdlgen1}). 
Since $B_1$ is analytic, we have
\[
DB_1(x-y)a^{(2)}\nu(y)=\sum_{j,r,z=1}^n\int_0^1 \frac{\partial B_1}{\partial x_j\partial x_r}
(s(x-y))\,ds(x_r-y_r)
a_{jz}\nu_z(y)
+DB_1(0)a^{(2)}\nu(y)
\]
for all $ (x,y)\in (\partial\Omega\times \partial\Omega)\setminus {\mathbb{D}}_{\partial\Omega} $ and 
\[
b \equiv  
\sup_{(x,y)\in (\partial\Omega)^2\setminus {\mathbb{D}}_{\partial\Omega}}
 \sum_{j,r=1}^n  \int_0^1 \left|\frac{\partial B_1}{\partial x_j\partial x_r}
(s(x-y)) \right| \,ds <+\infty\,.
\]
Hence,
\begin{eqnarray}\label{thm:mftgdlgech13}
 &&
 \biggl|\int_{\partial\Omega \setminus{\mathbb{B}}_n(x,r)  }
DB_1(x-y)a^{(2)}\nu(y) 
 \\ \nonumber
&&\quad 
\times\sum_{l=1}^n\nu_l(x)
\left[
\nu_l(x)\frac{x_h-y_h}{|x-y|^2}-\nu_h(x)\frac{x_l-y_l}{|x-y|^2}
\right]\,d\sigma_y\biggr|
\\ \nonumber
&& 
\leq 
 b
 |a^{(2)}|2nm_{n-1}(\partial\Omega)
 \\ \nonumber
&&\quad 
+|DB_1(0)a^{(2)}|\int_{\partial\Omega}|\nu(y)-\nu(x)|2n|x-y|^{-1}\,d\sigma_y 
\\ \nonumber
&&\quad 
+|DB_1(0)a^{(2)}\nu(x)|2\sum_{l=1}^n
\biggl|\int_{\partial\Omega \setminus{\mathbb{B}}_n(x,r)  }
\frac{x_l-y_l}{|x-y|^2}\,d\sigma_y\biggr|
\\ \nonumber
&&  
\leq 
 b
 |a^{(2)}|2nm_{n-1}(\partial\Omega)
+|DB_1(0)a^{(2)}|m_{n-1}(\partial\Omega){\mathrm{Lip}} (\nu) 2n 
\\ \nonumber
&&\quad 
+|DB_1(0)a^{(2)}|2n\sup_{l\in\{1,\dots,n\}}\sup_{\xi\in\partial\Omega}\sup_{\rho\in]0,+\infty[}
\biggl|\int_{\partial\Omega \setminus{\mathbb{B}}_n(\xi,\rho)  }
\frac{\xi_l-y_l}{|\xi-y|^2}\,d\sigma_y\biggr| 
\end{eqnarray}
for all $x\in \partial\Omega$ and $r\in]0,+\infty[$.  Then we note that if $n>2$, then 
\begin{eqnarray}\label{thm:mftgdlgech13a}
\lefteqn{\sup_{l\in\{1,\dots,n\}}\sup_{\xi\in\partial\Omega}\sup_{\rho\in]0,+\infty[}
\biggl|\int_{\partial\Omega \setminus{\mathbb{B}}_n(\xi,\rho)  }
\frac{\xi_l-y_l}{|\xi-y|^2}\,d\sigma_y\biggr|
}
\\ \nonumber
&&\quad 
\leq
\sup_{l\in\{1,\dots,n\}}\sup_{\xi\in\partial\Omega}\int_{\partial\Omega}
\left|
\frac{\xi_l-y_l}{|\xi-y|^2}
\right|\,d\sigma_y<+\infty 
\end{eqnarray}
(cf.~(\ref{thm:maggnfs2})).  If instead $n=2$, then we note that 
\[
 \frac{\xi_l}{|\xi|^2}\qquad\forall \xi\in 
{\mathbb{R}}^n\setminus\{0\}
\]
is odd and positively homogeneous of degree $-1=-(2-1)$ for all $l\in\{1,\dots,n\}$ and thus  Theorem \ref{lem:maolg} on the estimate of the maximal function associated to an odd positively homogeneous kernel of degree $-(2-1)$ implies that 
\begin{equation}\label{thm:mftgdlgech13b}
\sup_{l\in\{1,\dots,n\}}\sup_{\xi\in\partial\Omega}\sup_{\rho\in]0,+\infty[}
\biggl|\int_{\partial\Omega \setminus{\mathbb{B}}_n(\xi,\rho)  }
\frac{\xi_l-y_l}{|\xi-y|^2}\,d\sigma_y\biggr|<+\infty\,.
\end{equation}
  We now look at the eighth addendum $J_{8h}(x,r)$ in the right hand side of formula (\ref{lem:tgdlgen1}). Since the eighth addendum is bounded in $(\partial\Omega)^2\setminus {\mathbb{D}}_{\partial\Omega}$, we have
\begin{eqnarray}\label{thm:mftgdlgech14}
 &&
 \biggl|\int_{\partial\Omega \setminus{\mathbb{B}}_n(x,r)  }
\sum_{j,s=1}^n
\sum_{l=1}^n\nu_l(x)
\bigl[\nu_l(x)\frac{\partial^2C}{\partial x_h\partial x_j}(x-y)
\\ \nonumber
&&\quad 
-
\nu_h(x)\frac{\partial^2C}{\partial x_l\partial x_j}(x-y)
\bigr]
a_{js}\nu_s(y) \,d\sigma_y\biggr|
\\ \nonumber
&& 
\leq 2n^3\sup_{(x,y)\in(\partial\Omega)^2\setminus {\mathbb{D}}_{\partial\Omega}}
\sup_{h,j\in\{1,\dots,n\}}\biggl|\frac{\partial^2C}{\partial x_h\partial x_j}(x-y)\biggr|
|a^{(2)}|m_{n-1}(\partial\Omega)<+\infty
\end{eqnarray}
for all $x\in \partial\Omega$ and $r\in]0,+\infty[$.  We now look at the nineth addendum $J_{9h}(x,r)$ in the right hand side of formula (\ref{lem:tgdlgen1}). By formula (\ref{grafun1}), we have 
\[
C_{1, S_{ {\mathbf{a}} },\partial\Omega}\equiv
\sup_{ 0<|x|\leq{\mathrm{diam}}\,(\partial\Omega)}|\xi|^{n-1}| DS_{ {\mathbf{a}} } (\xi)|<+\infty 
\]
(see also \cite[Lem.~4.3 (ii)]{DoLa17}). By Theorem  \ref{thm:maggnfs}, we have
\begin{eqnarray} \label{thm:mftgdlgech15}
 &&
 \biggl|\int_{\partial\Omega \setminus{\mathbb{B}}_n(x,r)  }
\nu(y)^t\cdot a^{(1)} 
\sum_{l=1}^n\nu_l(x)
\bigl[\nu_l(x)\frac{\partial S_{{\mathbf{a}}} }{\partial x_h}(x-y)
- \nu_h(x)\frac{\partial S_{{\mathbf{a}}} }{\partial x_l}(x-y)
\bigr] \,d\sigma_y
\biggr|
\\  \nonumber 
&& 
\leq\biggl|\int_{\partial\Omega \setminus{\mathbb{B}}_n(x,r)  }
(\nu(y)-\nu(x))^t\cdot a^{(1)} 
\sum_{l=1}^n\nu_l(x)
\bigl[\nu_l(x)\frac{\partial S_{{\mathbf{a}}} }{\partial x_h}(x-y)
\\ \nonumber
&& 
- \nu_h(x)\frac{\partial S_{{\mathbf{a}}} }{\partial x_l}(x-y)
\bigr] \,d\sigma_y
\biggr|
\\  \nonumber 
&& 
+
\biggl|\int_{\partial\Omega \setminus{\mathbb{B}}_n(x,r)  }
\nu(x)^t\cdot a^{(1)} 
\sum_{l=1}^n\nu_l(x)
\bigl[\nu_l(x)\frac{\partial S_{{\mathbf{a}}} }{\partial x_h}(x-y)
- \nu_h(x)\frac{\partial S_{{\mathbf{a}}} }{\partial x_l}(x-y)
\bigr] \,d\sigma_y
\biggr|
\\  \nonumber 
&& 
\leq {\mathrm{Lip}}(\nu)|a^{(1)} |2n\sup_{l\in\{1,\dots,n\}}\sup_{\xi\in\partial\Omega}
\int_{\partial\Omega}|\xi-y|\left|\frac{\partial S_{{\mathbf{a}}} }{\partial x_l}(\xi-y)\right|\,d\sigma_y
\\  \nonumber 
&& 
+|a^{(1)} |2n\sup_{l\in\{1,\dots,n\}}
\sup_{\xi\in\partial\Omega}\sup_{\rho\in]0,+\infty[}
\biggl|\int_{\partial\Omega \setminus{\mathbb{B}}_n(\xi,\rho)  }
\frac{\partial S_{{\mathbf{a}}} }{\partial x_l}(\xi-y)\,d\sigma_y\biggr|<+\infty
\end{eqnarray}
for all $x\in \partial\Omega$ and $r\in]0,+\infty[$ (cf.~(\ref{thm:maggnfs2})). Hence, the proof of inequalities (\ref{thm:mftgdlgech2c}) is complete and accordingly (\ref{thm:dllreggenn1}) is equivalent to (\ref{thm:mftgdlgx1}).

\hfill  $\Box$ 

\vspace{\baselineskip}
 
 We note that by equality (\ref{eq:ppgr}), the kernel in parentheses of (\ref{thm:mftgdlgx1}) coincides with
the kernel $-\frac{1}{ \sqrt{\det a^{(2)} }}D(S_n\circ T^{-1})(x-y)a^{(2)}\nu (y)$, and that if the operator $P[D,{\mathbf{a}}]$ is homogenous, such a kernel is  precisely the kernel of the double layer potential. Next we prove   that condition
  (\ref{thm:dllreggenn1}) holds true.
\begin{theorem}\label{thm:mftgdlgx}
Let ${\mathbf{a}}$ be as in (\ref{introd0}), (\ref{ellip}), (\ref{symr}).  
Let $T\in M_{n}({\mathbb{R}})$ be an invertible matrix  as in (\ref{prop:ourfs0}).  Let $S_{ {\mathbf{a}} }$ be a fundamental solution of $P[{\mathbf{a}},D]$.  
Let $\Omega$ be a bounded open subset of class $C^{1,1}$ of ${\mathbb{R}}^n$.  Then  
condition (\ref{thm:dllreggenn1}) is satisfied.
\end{theorem}
  {\bf Proof.}    By Proposition \ref{prop:mftgdlgech}, it suffices to show  the validity of condition (\ref{thm:mftgdlgx1}). By the Lemma of the uniform cylinders, there exist   $r_{\partial\Omega}$, $\delta\in]0,1[$ with $r_{\partial\Omega}<\delta$ such that for each $\xi\in \partial\Omega$ there exists
$R_\xi\in O_n({\mathbb{R}})$ such that $C(\xi,R_\xi,r_{\partial\Omega},\delta)$ is a coordinate cylinder for $\partial\Omega$ around $\xi$ and that the corresponding function $\gamma_\xi$ satisfies  the conditions
\begin{equation}\label{thm:mftgdlgx0}
\gamma_\xi(0)=0\,,\quad
D\gamma_\xi(0)=0\qquad \forall\xi\in \partial\Omega\,,\quad
A\equiv\sup_{\xi\in \partial\Omega}\|\gamma_\xi\|_{ C^{1,1}(\overline{{\mathbb{B}}_{n-1}(0,r_{\partial\Omega})}) }<+\infty  
\end{equation}
(cf.~\textit{e.g.},   \cite[Lem.~2.63]{DaLaMu21} with Dalla Riva and Musolino).  Since   $\partial\Omega$ is compact and   the kernel  in the integral of (\ref{thm:mftgdlgx1})
 is continuous on the compact set $\{(x,y)\in (\partial\Omega)^2:\,|x-y|\geq r_{\partial\Omega}/2\}$ and an orthogonal map is an isometry,  condition (\ref{thm:mftgdlgx1}) is equivalent to
 \begin{eqnarray}\label{thm:mftgdlgx16}
&&\sup_{x\in \partial\Omega}\sup_{r\in ]0,r_{\partial\Omega}[}
\biggl\vert
\int_{[(\partial\Omega)\cap C(x,R_x,r_{\partial\Omega},\delta)]\setminus  {\mathbb{B}}_n(x,r)}
\\ \nonumber
&&\quad 
R_x{\mathrm{grad}}_{\partial\Omega,x}\left(
-\frac{1}{s_{n}\sqrt{\det a^{(2)} }}|T^{-1}(x-y)|^{-n} (x-y)^{t}\nu (y)
\right)\,d\sigma_y
\biggr\vert <+\infty\,,
\end{eqnarray}
an inequality that we now turn to verify. 
By formulas (\ref{eq:boest3}), (\ref{lem:tgdlgen1}), we have
\begin{eqnarray*} 
\lefteqn{
 \left( {\mathrm{grad}}_{\partial\Omega,x}\left(
-\frac{1}{s_{n}\sqrt{\det a^{(2)} }}|T^{-1}(x-y)|^{-n} (x-y)^{t}\nu (y)
\right) \right)_h
}
\\ \nonumber
&&\qquad\quad 
=\frac{n}{s_n\sqrt{\det a^{(2)} }}
 \frac{(x-y)^t\cdot\nu(y)}{|T^{-1}(x-y)|^n}
 \\ \nonumber
&&\qquad\quad 
 \times\sum_{l=1}^n\nu_l(x)
 \biggl[
 \nu_l(x)
 \frac{
 \sum_{j,z=1}^n (T^{-1})_{jz}(x_z-y_z)(T^{-1})_{jh}
  }{|T^{-1}(x-y)|^2}  
  \\ \nonumber
&&\qquad \quad 
 -
 \nu_h(x)\frac{\sum_{j,z=1}^n (T^{-1})_{jz}(x_z-y_z)(T^{-1})_{jl} 
  }{|T^{-1}(x-y)|^2}\biggr]
\\ \nonumber
&&\qquad\quad  
 -
 \frac{\sum_{l=1}^n\nu_l(x)
 \bigl[
 \nu_l(x) \nu_h(y)  
 -  \nu_h(x)  \nu_l(y) 
 \bigr]}{
 s_n\sqrt{\det a^{(2)} }   |T^{-1}(x-y)|^{n}} 
  \end{eqnarray*}
 for all $(x,y)\in (\partial\Omega)^2\setminus {\mathbb{D}}_{\partial\Omega}$
and $h\in\{1,\dots,n\}$.   We now wish to change the variables in the integral of (\ref{thm:mftgdlgx16}) with the map
\[
\psi_{x}(\eta)\equiv x+R_x^{t}\left(
\begin{array}{c}
\eta
\\
\gamma_x(\eta)
\end{array}
\right)\qquad\forall \eta\in {\mathbb{B}}_{n-1}(0,r_{\partial\Omega})\,.
\]
To do so, we note that
\begin{eqnarray*}
&&\nu(\psi_{x}(\eta) )=
R_x^t\left(
\begin{array}{c}
 -D\gamma_x(\eta)
\\
 1
\end{array}
\right) 
\frac{1}{\sqrt{1+\vert D\gamma_x(\eta)\vert^{2}}}
\quad\forall \eta\in {\mathbb{B}}_{n-1}(0,r_{\partial\Omega})\,, 
\\
&&
\nu(x)=R_x^t\left(
\begin{array}{c}
 0
\\
 1
\end{array}
\right) 
\end{eqnarray*}
for all $x\in\partial\Omega$ (cf.  \textit{e.g.}, \cite[Prop.~2.42]{DaLaMu21}) and that $(a^{(2)})^{-1}=(T^{-1})^tT^{-1}$. Then we have
\begin{eqnarray}\label{thm:mftgdlgx17}
&&
\int_{[(\partial\Omega)\cap C(x,R_x,r_{\partial\Omega},\delta)]\setminus  {\mathbb{B}}_n(x,r)}
\\ \nonumber
&&\quad 
\left(R_x{\mathrm{grad}}_{\partial\Omega,x}\left(
-\frac{1}{s_{n}\sqrt{\det a^{(2)} }}|T^{-1}(x-y)|^{-n} (x-y)^{t}\nu (y)
\right)\right)_b\,d\sigma_y
\\ \nonumber
&& 
=\int_{
\{\eta\in 
{\mathbb{B}}_{n-1}(0,r_{\partial\Omega}):\,
|\eta|^2+|\gamma_x(\eta)|^2\geq r^2\}
}
\frac{n}{s_n\sqrt{\det a^{(2)} }}
 \frac{(-\eta,-\gamma_x(\eta)) \left(
\begin{array}{c}
 -D\gamma_x(\eta)
\\
 1
\end{array}
\right) }{|T^{-1}R_x^t(-\eta,-\gamma_x(\eta))^t|^n}
 \\ \nonumber
&&\   
 \times\sum_{l,h=1}^n\left[R_x^t\left(
\begin{array}{c}
 0
\\
 1
\end{array}
\right) \right]_l 
 \biggl\{
 \left[R_x^t\left(
\begin{array}{c}
 0
\\
 1
\end{array}
\right) \right]_l 
 \frac{
 \sum_{j,s,z=1}^n (T^{-1})_{jz}
  (R_x^t)_{zs}(-\eta,-\gamma_x(\eta))^t_s
(T^{-1})_{jh}
  }{|T^{-1}R_x^t(-\eta,-\gamma_x(\eta))^t|^2}  
  \\ \nonumber
&& \ 
 -
\left[ R_x^t\left(
\begin{array}{c}
 0
\\
 1
\end{array}
\right) \right]_h \frac{\sum_{j,s,z=1}^n  (T^{-1})_{jz} (R_x^t)_{zs}(-\eta,-\gamma_x(\eta))^t_s(T^{-1})_{jl} 
  }{|T^{-1}R_x^t(-\eta,-\gamma_x(\eta))^t|^2}\biggr\} (R_x)_{bh}
\\ \nonumber
&& \   
 -
 \frac{\sum_{h=1}^n  (R_x)_{bh} \left[ R_x^t\left(
\begin{array}{c}
 -D\gamma_x(\eta)
\\
 1
\end{array}
\right) \right]_h 
 }{
 s_n\sqrt{\det a^{(2)} }   |T^{-1}R_x^t(-\eta,-\gamma_x(\eta))^t|^{n}}  
 \\ \nonumber
&& \   
 -
 \frac{
 -  \sum_{l,h=1}^n\left[ R_x^t\left(
\begin{array}{c}
 0
\\
 1
\end{array}
\right) \right]_l\left[ R_x^t\left(
\begin{array}{c}
 0
\\
 1
\end{array}
\right) \right]_h  \left[ R_x^t\left(
\begin{array}{c}
 -D\gamma_x(\eta)
\\
 1
\end{array}
\right) \right]_l 
( R_x)_{bh}}{
 s_n\sqrt{\det a^{(2)} }   |T^{-1}R_x^t(-\eta,-\gamma_x(\eta))^t|^{n}} \,d\eta
\\ \nonumber
&& 
=\int_{
\{\eta\in 
{\mathbb{B}}_{n-1}(0,r_{\partial\Omega}):\,
|\eta|^2+|\gamma_x(\eta)|^2\geq r^2\}
}
\frac{n}{s_n\sqrt{\det a^{(2)} }}
 \frac{D\gamma_x(\eta)\eta-\gamma_x(\eta)
 }{|T^{-1}R_x^t(\eta,\gamma_x(\eta))^t|^n}
 \\ \nonumber
&& \  
 \biggl[
 \frac{
 \sum_{h,s,z=1}^n 
  (R_x)_{bh}((a^{(2)})^{-1})_{hz}
  (R_x)^t_{zs}(-\eta,-\gamma_x(\eta))^t_s 
 }{|T^{-1}R_x^t(-\eta,-\gamma_x(\eta))^t|^2}  
  \\ \nonumber
&&  \  
 -
 \frac{\sum_{l,s,z=1}^n 
\left(
\begin{array}{c}
 0
\\
 1
\end{array}
\right)_b (R_x)_{nl}((a^{(2)})^{-1})_{lz}
 (R_x)^t_{zs}(-\eta,-\gamma_x(\eta))^t_s 
  }{|T^{-1}R_x^t(-\eta,-\gamma_x(\eta))^t|^2} \biggr] 
  \\ \nonumber
&& \   
 -
 \frac{ 
    \left(
\begin{array}{c}
 -D\gamma_x(\eta)
\\
 1
\end{array}
\right)_b  
 -     \left(
\begin{array}{c}
 0
\\
 1
\end{array}
\right)_b
   }{
 s_n\sqrt{\det a^{(2)} }   |T^{-1}R_x^t(-\eta,-\gamma_x(\eta))^t|^{n}} \,d\eta\qquad\forall r\in]0,r_{\partial\Omega}[ 
\end{eqnarray}
for all $b\in\{1,\dots,n\}$. Hence,
\begin{eqnarray}\label{thm:mftgdlgx18}
&&
\int_{[(\partial\Omega)\cap C(x,R_x,r_{\partial\Omega},\delta)]\setminus  {\mathbb{B}}_n(x,r)}
\\ \nonumber
&&\qquad 
\left(R_x{\mathrm{grad}}_{\partial\Omega,x}\left(
-\frac{1}{s_{n}\sqrt{\det a^{(2)} }}|T^{-1}(x-y)|^{-n} (x-y)^{t}\nu (y)
\right)\right)_b\,d\sigma_y
\\ \nonumber
&&\quad
=0\qquad \forall r\in]0,r_{\partial\Omega}[ 
\end{eqnarray}
if $b=n$,
\begin{eqnarray}\label{thm:mftgdlgx18a}
&&
\int_{[(\partial\Omega)\cap C(x,R_x,r_{\partial\Omega},\delta)]\setminus  {\mathbb{B}}_n(x,r)}
\\ \nonumber
&&\qquad 
\left(R_x{\mathrm{grad}}_{\partial\Omega,x}\left(
-\frac{1}{s_{n}\sqrt{\det a^{(2)} }}|T^{-1}(x-y)|^{-n} (x-y)^{t}\nu (y)
\right)\right)_b\,d\sigma_y
\\ \nonumber
&&\quad
=\int_{
\{\eta\in 
{\mathbb{B}}_{n-1}(0,r_{\partial\Omega}):\,
|\eta|^2+|\gamma_x(\eta)|^2\geq r^2\} 
}
F_b(x,\eta)\,d\eta\qquad  \forall r\in]0,r_{\partial\Omega}[ 
\end{eqnarray}
if $ b\in\{1,\dots,n-1\}$, where
\begin{eqnarray*}
\lefteqn{
F_b(x,\eta)
\equiv
\frac{n}{s_n\sqrt{\det a^{(2)} }}
 \frac{D\gamma_x(\eta)\eta-\gamma_x(\eta)
 }{|T^{-1}R_x^t(\eta,\gamma_x(\eta))^t|^n}
 }
 \\ \nonumber
&&\qquad\quad 
 \times\biggl[
 \frac{
 \sum_{h,s,z=1}^n 
 (R_x)_{bh}((a^{(2)})^{-1})_{hz}
 (R_x^t)_{zs}(-\eta,-\gamma_x(\eta))^t_s 
 }{|T^{-1}R_x^t(-\eta,-\gamma_x(\eta))^t|^2}  \biggr] 
  \\ \nonumber
&&\qquad\quad  
 +
 \frac{ 
 \frac{\partial\gamma_x}{\partial\eta_b}(\eta)
}{
 s_n\sqrt{\det a^{(2)} }   |T^{-1}R_x^t(-\eta,-\gamma_x(\eta))^t|^{n}}  \qquad 
 \forall \eta\in {\mathbb{B}}_{n-1}(0,r_{\partial\Omega})
\end{eqnarray*}
for all $x\in\partial\Omega$ and $b\in\{1,\dots,n-1\}$. By assumption, we have $\gamma_x(0)=0$, $D\gamma_x(0)=0$ and thus the Taylor formula with integral residue implies that
\begin{equation}\label{thm:mftgdlgx19}
|\gamma_x(\eta)|\leq
 | D \gamma_x |_{1}
   |\eta|^{2}
   \leq A |\eta|^{2}\,\quad
   | D \gamma_x (\eta)|\leq | D \gamma_x |_{1}|\eta|\leq A|\eta|
\end{equation}
for all $\eta\in{\mathbb{B}}_{n-1}(0,r_{\partial\Omega})$ and $h\in\{1,\dots,n-1\}$  (cf.~\textit{e.g.}, \cite[Lem.~2.58]{DaLaMu21}). Then we set
\[
a_x(r)\equiv\left\{
\frac{|\gamma_x(\eta)|}{|\eta|}:\,\eta\in \overline{{\mathbb{B}}_{n-1}(0,r)}\setminus\{0\}
\right\}\qquad\forall r\in]0,r_{\partial\Omega}[ 
\]
and we note that
\[
a_x(r)\leq | D \gamma_x |_{1} r 
\leq A  r 
\qquad\forall r\in]0,r_{\partial\Omega}[ 
\]
for all $x\in\partial\Omega$. Then the triangular inequality, inequality (\ref{thm:bmfdl3}), inequality $|R_x|\leq 1$  and inequality 
\[
|T^{-1}R_x^t(\eta,\gamma_x(\eta))^t|\geq |T|^{-1} |(\eta,\gamma_x(\eta))|\geq |T|^{-1}  |\eta|\qquad\forall\eta\in {\mathbb{B}}_{n-1}(0,r_{\partial\Omega})
\]
for all $x\in\partial\Omega$ imply that
\[
c_F\equiv\sup_{b\in\{1,\dots,n-1\}}\sup_{x\in\partial\Omega}\sup_{\eta\in {\mathbb{B}}_{n-1}(0,r_{\partial\Omega})\setminus\{0\}}\frac{|F_b(x,\eta)|}{|\eta|^{n-1}}
<+\infty \,.
\]
 Then Lemma \ref{lem:pvcom} (iii) of the Appendix implies that
\begin{eqnarray}\label{thm:mftgdlgx20}
\lefteqn{ 
\sup_{b\in\{1,\dots,n-1\}}\sup_{x\in\partial\Omega}\biggl| \int_{
 \{\eta\in 
{\mathbb{B}}_{n-1}(0,r_{\partial\Omega}):\,
|\eta|^2+|\gamma_x(\eta)|^2\geq r^2\} 
}
F_b(x,\eta)\,d\eta
}
\\ \nonumber
&&\qquad\quad
-
 \int_{
 \{\eta\in 
{\mathbb{B}}_{n-1}(0,r_{\partial\Omega}):\,
|\eta|^2\geq r^2\} 
}
F_b(x,\eta)\,d\eta\biggr|
\\ \nonumber
&&
\leq c_Fs_{n-1}
\log
(1+\alpha(r))
\leq c_Fs_{n-1}
\log (1+A r_{\partial\Omega})
\qquad\forall r\in]0,r_{\partial\Omega}[\,.
\end{eqnarray}
Next  note that 
\[
F_b(x,\eta)=\tilde{F}_b(x,\eta)+F^\sharp_b(x,\eta) \qquad   \forall \eta\in {\mathbb{B}}_{n-1}(0,r_{\partial\Omega})\,,
\]
where
\begin{eqnarray*} 
\lefteqn{
\tilde{F}_b(x,\eta)\equiv\frac{n}{s_n\sqrt{\det a^{(2)} }}
 \frac{D\gamma_x(\eta)\eta-2\gamma_x(\eta)
 }{|T^{-1}R_x^t(\eta,\gamma_x(\eta))^t|^n}
 }
 \\ \nonumber
&&\qquad\quad 
\times \biggl[
 \frac{
 \sum_{h,s,z=1}^n 
 (R_x)_{bh}((a^{(2)})^{-1})_{hz}
 (R_x)^t_{zs}(-\eta,-\gamma_x(\eta))^t_s 
 }{|T^{-1}R_x^t(-\eta,-\gamma_x(\eta))^t|^2}  \biggr]\,,
 \\ \nonumber
\lefteqn{
F^\sharp_b(x,\eta)\equiv\frac{n}{s_n\sqrt{\det a^{(2)} }}
 \frac{\gamma_x(\eta)
 }{|T^{-1}R_x^t(\eta,\gamma_x(\eta))^t|^n}
 }
 \\ \nonumber
&&\qquad\quad 
\times \biggl[
 \frac{
 \sum_{h,s,z=1}^n 
 (R_x)_{bh}((a^{(2)})^{-1})_{hz}
 (R^t_x)_{zs}(-\eta,-\gamma_x(\eta))^t_s 
 }{|T^{-1}R_x^t(-\eta,-\gamma_x(\eta))^t|^2}  \biggr]  
\\ \nonumber
&&\qquad\quad 
+
\frac{ 
 \frac{\partial\gamma_x}{\partial\eta_b}(\eta)
}{
 s_n\sqrt{\det a^{(2)} }   |T^{-1}R_x^t(-\eta,-\gamma_x(\eta))^t|^{n}}  \qquad   \forall \eta\in {\mathbb{B}}_{n-1}(0,r_{\partial\Omega})
\end{eqnarray*}
for all $x\in\partial\Omega$ and $b\in\{1,\dots,n-1\}$. Next we note that the function
\begin{eqnarray*}
\lefteqn{
k_{x,b}(y)\equiv\frac{n}{s_n\sqrt{\det a^{(2)} }}
 \frac{1
 }{|T^{-1}R_x^t(y)^t|^n}
 }
 \\ \nonumber
&&\qquad\quad 
 \times	 
 \biggl[
 \frac{
 \sum_{h,s,z=1}^n 
 (R_x)_{bh}((a^{(2)})^{-1})_{hz}
 (R^t_x)_{zs}(-y)_s 
 }{|T^{-1}R_x^ty|^2}  \biggr]\qquad\forall y\in {\mathbb{R}}^n\setminus\{0\}
 \end{eqnarray*}
is positively homogeneous of degree $-(n+1)$ and odd. Moreover,  the triangular inequality, inequality (\ref{thm:bmfdl3}), inequality $|R_x|\leq 1$  and inequality 
\[
|T^{-1}R_x^ty|\geq |T|^{-1}|R_x^ty|= |T|^{-1} |y|\qquad\forall y\in {\mathbb{R}}^n\setminus\{0\}
\]
for all $x\in\partial\Omega$ imply that
\[
\sup_{b\in\{1,\dots,n-1\}}\sup_{x\in\partial\Omega}\sup_{
y\in {\mathbb{B}}_{n}(0,3/2)\setminus {\mathbb{B}}_{n}(0,1/2)
}|D_yk_{x,b}(y)|<+\infty\,.
\]
 Since the restriction from $C^1(\overline{{\mathbb{B}}_{n}(0,3/2)}\setminus {\mathbb{B}}_{n}(0,1/2))$
 to $C^{0,1}(\partial{\mathbb{B}}_{n}(0,1))$ is linear and continuous, we have
 \[
 K\equiv\sup_{b\in\{1,\dots,n-1\}}\sup_{x\in\partial\Omega}\| k_{x,b}\|_{C^{0,1}(\partial{\mathbb{B}}_{n}(0,1))}
 <+\infty\,.
 \]
 Then by a known inequality for positively homogeneous functions of degree $-(n+1) $ and by inequality (\ref{thm:mftgdlgx19}), we have
 \begin{eqnarray*}
\lefteqn{
 |k_{x,b}(\eta,\gamma_x(\eta))-k_{x,b}(\eta,0)|
}
\\ \nonumber
&&\qquad
\leq  ( 2^1+2(n+1))
 \max\{\sup_{ \partial{\mathbb{B}}_n(0,1) }|k_{x,b}|, |k_{x,b}:\,\partial{\mathbb{B}}_n(0,1)|_1\}
 \\ \nonumber
&&\qquad\quad
 \times |(\eta,\gamma_x(\eta))-(\eta,0)|^1(\min\{|(\eta,\gamma_x(\eta))|,|(\eta,0)|\})^{-(n+1)-1}
 \\ \nonumber
&&\qquad
\leq (2n+
 4	 
)KA|\eta|^2
   |\eta|^{-n-2} \qquad\forall \eta\in {\mathbb{B}}_{n-1}(0,r_{\partial\Omega})\setminus\{0\}
 \end{eqnarray*}
(see Cialdea \cite[VIII, p.~47]{Ci95}, or also 
\cite[Lem.~4.14]{DaLaMu21}, case $\alpha=1$). Then again by  (\ref{thm:mftgdlgx19}), we have
\[
|D\gamma_x(\eta)\eta-2\gamma_x(\eta)|\leq 3A|\eta|^2
 \qquad\forall\eta\in {\mathbb{B}}_{n-1}(0,r_{\partial\Omega})\setminus\{0\}\,.
\]
Hence,
 \begin{eqnarray}\label{thm:mftgdlgx22}
&&
\biggl|
\int_{
 \{\eta\in 
{\mathbb{B}}_{n-1}(0,r_{\partial\Omega}):\,
|\eta|^2\geq r^2\} 
}
\tilde{F}_b(x,\eta)\,d\eta  \biggr|
\\ \nonumber
&& 
\leq \biggl|
\int_{
 \{\eta\in 
{\mathbb{B}}_{n-1}(0,r_{\partial\Omega}):\,
|\eta|^2\geq r^2\} 
}
\tilde{F}_b(x,\eta)   
 -(D\gamma_x(\eta)\eta-2\gamma_x(\eta))k_{x,b}(\eta,0)\,d\eta\biggr|
 \\ \nonumber
&& \qquad\qquad 
 +\biggl|\int_{
 \{\eta\in 
{\mathbb{B}}_{n-1}(0,r_{\partial\Omega}):\,
|\eta|^2\geq r^2\} 
}(D\gamma_x(\eta)\eta-2\gamma_x(\eta))k_{x,b}(\eta,0)\,d\eta\biggr|
 \\ \nonumber
&&
 \leq  
\int_{
{\mathbb{B}}_{n-1}(0,r_{\partial\Omega})
}
3A|\eta|^2
  (2n+ 4	 )KA|\eta|^{-n}  \,d\eta 
 \\ \nonumber
&& \qquad\qquad 
 +\biggl|\int_{
 \{\eta\in 
{\mathbb{B}}_{n-1}(0,r_{\partial\Omega}):\,
|\eta|^2\geq r^2\} 
} (D\gamma_x(\eta)\eta-2\gamma_x(\eta))k_{x,b}(\eta,0)\,d\eta\biggr|
\end{eqnarray}
for all $ r\in]0,r_{\partial\Omega}[$.
 Next we estimate the last integral in the right hand side and we distinguish cases $n>2$ and $n=2$. If $n>2$, we integrate on the spheres and we obtain  again by  (\ref{thm:mftgdlgx19})
 \begin{eqnarray}\label{thm:mftgdlgx23}
\lefteqn{
\biggl|\int_{
 \{\eta\in 
{\mathbb{B}}_{n-1}(0,r_{\partial\Omega}):\,
|\eta|^2\geq r^2\} 
} (D\gamma_x(\eta)\eta-2\gamma_x(\eta))k_{x,b}(\eta,0)\,d\eta\biggr|
}
\\ \nonumber
&& 
=\biggl|\int_{ 
\partial{\mathbb{B}}_{
 n-1	 
}(0,1) 
} \int_r^{r_{\partial\Omega}}
(D\gamma_x(\tau v)\tau v-2\gamma_x(\tau v))k_{x,b}(\tau v,0) \tau^{n-2}\,d\tau \,d\sigma_v\biggr|
\\ \nonumber
&& 
=\biggl|\int_{ 
\partial{\mathbb{B}}_{
 n-1	 
}(0,1) 
} \int_r^{r_{\partial\Omega}}
\frac{d}{d\tau}\left(\frac{\gamma_x(\tau v)}{\tau^2}\right)\tau^3\tau^{-(n+1)}k_{x,b}( v,0) \tau^{n-2} \,d\tau\,d\sigma_v\biggr|
\\ \nonumber
&& 
=\biggl|\int_{ 
\partial{\mathbb{B}}_{
 n-1	 
}(0,1) 
} \biggl[
 \frac{\gamma_x(\tau v)}{\tau^2} \biggr]_{\tau=r}^{\tau=r_{\partial\Omega}} k_{x,b}( v,0)  \,d\sigma_v\biggr|
 \\ \nonumber
&& 
 \leq \int_{ 
\partial{\mathbb{B}}_{
 n-1	 
}(0,1) }2A|v|^2K|v|^{-(n+1)}  \,d\sigma_v 
=2AKs_{
 n-1	 
}
\quad \forall    r\in]0,r_{\partial\Omega}[
\end{eqnarray}
for all $x\in\partial\Omega$ and $b\in\{1,\dots,n-1\}$.  If   instead $n=2$, we have
 \begin{eqnarray}\label{thm:mftgdlgx23a}
\lefteqn{
\biggl|\int_{
 \{\eta\in 
{\mathbb{B}}_{n-1}(0,r_{\partial\Omega}):\,
|\eta|^2\geq r^2\} 
} (D\gamma_x(\eta)\eta-2\gamma_x(\eta))k_{x,b}(\eta,0)\,d\eta\biggr|
}
\\ \nonumber
&&
=\biggl|\int_{]-r_{\partial\Omega},r_{\partial\Omega}[\setminus]-r,r[}
\frac{d}{d\eta}\left(\frac{\gamma_x(\eta)}{\eta^2}\right)\eta^3|\eta|^{-(2+1)}k_{x,b}( \eta/|\eta|,0) \,d\eta\biggr|
\\ \nonumber
&&\leq
\biggl|\biggl[\frac{\gamma_x(\eta)}{\eta^2}(-1)k_{x,b}(-1,0) 
\biggr]_{\eta=-r_{\partial\Omega}}^{\eta=-r}\biggr|
+\biggl|\biggl[\frac{\gamma_x(\eta)}{\eta^2}(1)k_{x,b}(1,0) 
\biggr]_{\eta=r}^{\eta=r_{\partial\Omega}}\biggr|
\\ \nonumber
&&
\leq 2AKs_{
 n-1	 
}\qquad  \forall    r\in]0,r_{\partial\Omega}[
\end{eqnarray}
for all $x\in\partial\Omega$ and $b=1$, where $s_{
 n-1	
}=s_{
2-
 1	 
}\equiv 2$. Next we turn to consider the integral of $F^\sharp_b(x,\eta)$ and we note that
\begin{eqnarray*}
\lefteqn{
\frac{1}{s_n\sqrt{\det a^{(2)} }}\frac{\partial}{\partial\eta_b}\left(
\frac{\gamma_x(\eta)}{|T^{-1}R_x^t(\eta,\gamma_x(\eta))^t|^n}\right)
}
\\ \nonumber
&&\qquad
=F^\sharp_b(x,\eta)
-
\frac{n}{s_n\sqrt{\det a^{(2)} }}
 \frac{\gamma_x(\eta)
 }{|T^{-1}R_x^t(\eta,\gamma_x(\eta))^t|^n}
\\ \nonumber
&&\qquad\quad 
\times \biggl[
 \frac{
 \sum_{h,s,z=1}^n  
 \frac{\partial\gamma_x}{\partial\eta_b}(\eta)
 (R_x)_{nh}((a^{(2)})^{-1})_{hz}
 (R^t_x)_{zs}(\eta,\gamma_x(\eta))^t_s 
 }{|T^{-1}R_x^t(\eta,\gamma_x(\eta))^t|^2}  \biggr]   
\end{eqnarray*}
for all $\eta \in {\mathbb{B}}_{n-1}(0,r_{\partial\Omega})\setminus\{0\} $,
$x\in\partial\Omega$ and $b\in\{1,\dots,n-1\}$.   If $n>2$, then the Divegence Theorem and inequalities (\ref{thm:bmfdl3}), (\ref{thm:mftgdlgx19})  imply that
\begin{eqnarray}\label{thm:mftgdlgx24}
&&
\biggl|\int_{
 \{\eta\in 
{\mathbb{B}}_{n-1}(0,r_{\partial\Omega}):\,
|\eta|^2\geq r^2\} 
}F^\sharp_b(x,\eta)\,d\eta\biggr|
\\ \nonumber
&& 
\leq\frac{1}{s_n\sqrt{\det a^{(2)} }}\biggl|\int_{\partial{\mathbb{B}}_{n-1}(0,r_{\partial\Omega})}\frac{\eta_b}{|\eta|}\frac{\gamma_x(\eta)\,d\sigma_\eta}{|T^{-1}R_x^t(\eta,\gamma_x(\eta))^t|^n}
\\ \nonumber
&&\qquad 
-
\int_{\partial{\mathbb{B}}_{n-1}(0,r )}\frac{\eta_b}{|\eta|}\frac{\gamma_x(\eta)\,d\sigma_\eta}{|T^{-1}R_x^t(\eta,\gamma_x(\eta))^t|^n}\biggr|
\\ \nonumber
&&\qquad 
+
\int_{{\mathbb{B}}_{n-1}(0,r_{\partial\Omega})}
\frac{n}{s_n\sqrt{\det a^{(2)} }}\frac{A|\eta|^2}{|T|^{-n}|\eta|^n}\frac{A|\eta||(a^{(2)})^{-1}|(|\eta|+A|\eta|^{
 2	 
})}{|T|^{-2}|\eta|^2}
\,d\eta
\\ \nonumber
&& 
\leq \frac{1}{s_n\sqrt{\det a^{(2)} }}2s_{
 n-1	 
}|T|^{n}A 
+\frac{
 n 	 
(A+1)^3|(a^{(2)})^{-1}||T|^{n+2}}{s_n\sqrt{\det a^{(2)} }}\int_{{\mathbb{B}}_{n-1}(0,r_{\partial\Omega})}
 \frac{d\eta}{|\eta|^{n-2}} 
\end{eqnarray}
for all   $x\in\partial\Omega$, $b\in\{1,\dots,n-1\}$ and  $r\in]0,r_{\partial\Omega}[$. Case $n=2$ can be treated similarly, with the understanding that $s_{
 n-1	 
}=2$ if $n=2$.  Then equalities (\ref{thm:mftgdlgx17}), (\ref{thm:mftgdlgx18}), (\ref{thm:mftgdlgx18a}), inequalities
(\ref{thm:mftgdlgx20}), (\ref{thm:mftgdlgx22}), (\ref{thm:mftgdlgx23}), (\ref{thm:mftgdlgx23a}) and (\ref{thm:mftgdlgx24})
imply the validity of condition (\ref{thm:mftgdlgx16}) and thus the proof is complete.\hfill  $\Box$ 

\vspace{\baselineskip}

\appendix

\section{
 A 	 
classical lemma for the computation of the principal value}\label{sec:app}
We now present the following variant of a known technical lemma that we have exploited in the paper. For the convenience of the reader, we include a proof, that is based on ideas of   Gegelia~\cite{Ge63}, Kupradze,  Gegelia,  Basheleishvili and 
 Burchuladze~\cite[p.~208]{KuGeBaBu79}, Mikhlin \cite[p.~40]{Mikh65}).
\begin{lemma}\label{lem:pvcom}
 Let $r\in]0,+\infty[$. Let $\gamma\in 
  C^0 
 ({\mathbb{B}}_{n-1}(0,r))$ 
 be bounded.	 
  Let $\gamma$ be differentiable at $0$. Let $\gamma(0)=0$. Let
\[
a\equiv (a_1,\dots,a_{n-1})\equiv D\gamma(0) \,. 
\]
Let
\[
q_{0,1}(\gamma)\equiv \sup\{|\eta|^{-1}|\gamma(\eta)|:\,
\eta\in {\mathbb{B}}_{n-1}(0,r)\setminus\{0\}
\}\,.
\]
Let 
\begin{eqnarray}\label{lem:pvcom1}
 A_a(\epsilon)&\equiv&\{\eta\in {\mathbb{B}}_{n-1}(0,r):\,|\eta|^2+|a\cdot\eta|^2<\epsilon^2\}\,, 
\\ \nonumber
A_{\gamma}(\epsilon)&\equiv&  \{\eta\in {\mathbb{B}}_{n-1}(0,r)\,:
|\eta|^2+|\gamma(\eta)|^2<\epsilon^2\}\qquad\forall \epsilon\in ]0,r[ \,,
 \end{eqnarray}
and
\[
 {\mathcal{A}}_{a} \equiv \{  A_{a}(\epsilon):\, \epsilon\in]0,r[\} \,,
 \qquad{\mathcal{A}}_{\gamma} \equiv \{  A_{\gamma}(\epsilon):\, \epsilon\in]0,r[\}\,.
\]
Then
\begin{enumerate}
\item[(i)] $q_{0,1}(\gamma)<+\infty$. 
\item[(ii)] 
\begin{eqnarray*}
 &&{\mathbb{B}}_{n-1}\left(0, \frac{\epsilon}{\sqrt{1+|a|^2}}\right)
\subseteq A_{a}(\epsilon)\subseteq {\mathbb{B}}_{n-1}(0,\epsilon)\,,
\\ \nonumber
&&
{\mathbb{B}}_{n-1}\left(0, \frac{\epsilon}{\sqrt{1+q_{0,1}(\gamma)^2}}\right)
\subseteq A_{\gamma }(\epsilon)\subseteq {\mathbb{B}}_{n-1}(0,\epsilon)
\end{eqnarray*}
for all $\epsilon\in ]0,r[$.
\item[(iii)] Let $g$ be a measurable function from ${\mathbb{B}}_{n-1}(0,r)$ to ${\mathbb{C}}$ such that $g$ is integrable in ${\mathbb{B}}_{n-1}(0,r)\setminus {\mathbb{B}}_{n-1}(0,\epsilon)$ for all $\epsilon\in]0,r[$, and such that there exists
$r_g\in]0,r]$ such that
\[
c_g\equiv{\mathrm{ess}}\,\sup\{|\eta|^{n-1}|g(\eta)|:\, \eta\in {\mathbb{B}}_{n-1}(0,r_g)\setminus\{0\}\}<+\infty\,.
\]
Let 
\begin{equation}\label{lem:pvcom1a}
\alpha(\epsilon)\equiv \sup\left\{
\frac{|\gamma(\eta)-a\cdot\eta|}{|\eta|}:\,\eta\in  {\mathbb{B}}_{n-1}(0,\epsilon)\setminus\{0\}
\right\}\qquad\forall \epsilon\in ]0,r[\,.
\end{equation}
Then 
\begin{eqnarray}
\label{lem:pvcom1b}
\lefteqn{
\left|
\int_{{\mathbb{B}}_{n-1}(0,r)\setminus A_\gamma(\epsilon)}g(\eta)\,d\eta
- \int_{{\mathbb{B}}_{n-1}(0,r)\setminus A_a(\epsilon)}g(\eta)\,d\eta
\right|
}
\\ \nonumber
&&\qquad
\leq c_gs_{n-1}
\log
(1+\alpha(\epsilon))
+c_gs_{n-1}\log\sqrt{1+|a|^2}
\end{eqnarray}
for all $\epsilon\in ]0, r_g [$  and
\begin{eqnarray}
\label{lem:pvcom1c}
\lefteqn{
\left|
\int_{{\mathbb{B}}_{n-1}(0,r)\setminus A_\gamma(\epsilon)}g(\eta)\,d\eta
- \int_{{\mathbb{B}}_{n-1}(0,r)\setminus A_a(\epsilon)}g(\eta)\,d\eta
\right|
}
\\ \nonumber
&&\qquad
\leq c_gs_{n-1}
\log
(1+\alpha(\epsilon))
+c_gs_{n-1}\log\frac{1}{ \sqrt{1-2\alpha(\epsilon)} }
\end{eqnarray}
for all $\epsilon\in ]0, r_g [$ such that $\alpha(\epsilon)<1/2$,  
with the understanding that $s_1=2$.
In particular, the limit
\[
\lim_{\epsilon\to 0}\int_{{\mathbb{B}}_{n-1}(0,r)\setminus A_\gamma(\epsilon)}g(\eta)\,d\eta
\]
exists if and only if the limit
\[
\lim_{\epsilon\to 0}\int_{{\mathbb{B}}_{n-1}(0,r)\setminus A_a(\epsilon)}g(\eta)\,d\eta
\]
exists, and in case of existence, the value of  two limits is equal. 
\end{enumerate}
\end{lemma}
{\bf Proof.} (i) By the differentiability of $\gamma $ at $0$, we have 
\begin{equation}\label{lem:pvcom2}
\lim_{\eta\to 0}\frac{|\gamma(\eta)-a\cdot\eta|}{|\eta|}=0\,,
\end{equation}
and thus there exists $r_1\in]0,r[$ such that
\[
|\gamma(\eta)-a\cdot\eta|\leq|\eta|\qquad\forall \eta\in {\mathbb{B}}_{n-1}(0,r_1)\setminus\{0\}\,.
\]
Hence, 
\[
|\gamma(\eta)|\leq |\eta|(|a|+1)\qquad\forall \eta\in {\mathbb{B}}_{n-1}(0,r_1)\setminus\{0\}\,.
\]
Since $\gamma$ and $|\eta|^{-1}$ are bounded on ${\mathbb{B}}_{n-1}(0,r  )\setminus {\mathbb{B}}_{n-1}(0,r_1 )$, we conclude that $q_{0,1}(\gamma)<+\infty$. 

(ii) If  
\[
\eta\in {\mathbb{B}}_{n-1}\left(0, \frac{\epsilon}{\sqrt{1+|a|^2}}\right)\,,
\]
then we have 
\[
|\eta|^2+|a\cdot\eta|^2\leq |\eta|^2+|a|^2|\eta|^2< \epsilon^2\,,
\]
and thus $\eta\in A_{a}(\epsilon)$. By   definition of $ A_{a}(\epsilon)$, we have
$  A_{a}(\epsilon)\subseteq {\mathbb{B}}_{n-1}(0,\epsilon)$.  If  
\[
\eta\in {\mathbb{B}}_{n-1}\left(0, \frac{\epsilon}{\sqrt{1+q_{0,1}(\gamma)^2}}\right)\,,
\]
then we have 
\[
|\eta|^2+|\gamma(\eta)|^2\leq |\eta|^2+q_{0,1}(\gamma)^2|\eta|^2< \epsilon^2\,,
\]
and thus $\eta\in A_{\gamma}(\epsilon)$. By  definition of $  A_{\gamma}(\epsilon)$, we have
$ A_{\gamma}(\epsilon)\subseteq {\mathbb{B}}_{n-1}(0,\epsilon)$. 

(iii) We first  estimate the difference
\begin{eqnarray}
\label{lem:pvcom3}
\lefteqn{
\left|
\int_{{\mathbb{B}}_{n-1}(0,r)\setminus A_\gamma(\epsilon)}g(\eta)\,d\eta
- \int_{{\mathbb{B}}_{n-1}(0,r)\setminus A_a(\epsilon)}g(\eta)\,d\eta
\right|
}
\\ \nonumber
&&\qquad
\leq
\int_{
\{ \eta\in  {\mathbb{B}}_{n-1}(0,r):\,   \epsilon^2\leq |\eta|^2+|\gamma(\eta)|^2\,,
\      |\eta|^2+|a\cdot\eta|^2<\epsilon^2    \}
}|g(\eta)|\,d\eta
\\ \nonumber
&&\qquad\quad
+
\int_{
\{ \eta\in  {\mathbb{B}}_{n-1}(0,r):\,   \epsilon^2\leq
|\eta|^2+|a\cdot\eta|^2\,,\ 
|\eta|^2+|\gamma(\eta)|^2  <\epsilon^2    \}
}|g(\eta)|\,d\eta 
\end{eqnarray}
 for each $\epsilon\in]0,r_g[$.	 
 In order to estimate first integral in the right hand side of (\ref{lem:pvcom3}), we prove the following inclusion
\begin{eqnarray}\label{lem:pvcom4}
&&
\{ \eta\in  {\mathbb{B}}_{n-1}(0,r):\,   \epsilon^2\leq |\eta|^2+|\gamma(\eta)|^2\,,
\      |\eta|^2+|a\cdot\eta|^2<\epsilon^2    \}
\\ \nonumber
&&\quad
\subseteq
\{ \eta\in  {\mathbb{B}}_{n-1}(0,r):\, 
\left(\frac{\epsilon}{1+\alpha(\epsilon)}\right)^2\leq   |\eta|^2+|a\cdot\eta|^2<\epsilon^2    \}
\qquad\forall\epsilon\in]0,r[\,.
\end{eqnarray}
To do so, we note that if $ \epsilon^2\leq |\eta|^2+|\gamma(\eta)|^2$ and $ |\eta|^2+|a\cdot\eta|^2<\epsilon^2$, then we have
\begin{eqnarray*}
\lefteqn{
 \epsilon^2\leq |\eta|^2+|\gamma(\eta)|^2
 =
|\eta|^2+(a\cdot\eta +(\gamma(\eta)-a\cdot\eta))^2
}
\\ \nonumber
&&\qquad
\leq  |\eta|^2+\left(  |a\cdot\eta|+\alpha(\epsilon)|\eta|  \right)^2
\\ \nonumber
&&\qquad
= |\eta|^2+ |a\cdot\eta|^2+\alpha(\epsilon)^2|\eta| ^2+2 |a\cdot\eta|\alpha(\epsilon)|\eta|
\\ \nonumber
&&\qquad
=(|\eta|^2+ |a\cdot\eta|^2)
\left(
1+\frac{\alpha(\epsilon)^2|\eta| ^2+2 |a\cdot\eta|\alpha(\epsilon)|\eta|   }{|\eta|^2+ |a\cdot\eta|^2}
\right)
\\ \nonumber
&&\qquad
\leq (|\eta|^2+ |a\cdot\eta|^2)\left(
1+\alpha(\epsilon)^2+2 \alpha(\epsilon)
\right)=(|\eta|^2+ |a\cdot\eta|^2)(1+\alpha(\epsilon))^2\,.
\end{eqnarray*}
In order to estimate the second integral in the right hand side of (\ref{lem:pvcom3})
we note that
\begin{eqnarray}\label{lem:pvcom5a}
\lefteqn{
\{ \eta\in  {\mathbb{B}}_{n-1}(0,r):\,   \epsilon^2\leq
|\eta|^2+|a\cdot\eta|^2\,,\ 
|\eta|^2+|\gamma(\eta)|^2  <\epsilon^2    \}
}
\\ \nonumber
&&\qquad
\subseteq
\{ \eta\in  {\mathbb{B}}_{n-1}(0,r):\,   \epsilon^2\leq
|\eta|^2+|a\cdot\eta|^2\,,\ |\eta|^2<\epsilon^2
\}
\qquad\forall \epsilon\in]0,r[\,.
\end{eqnarray} 
In order to estimate second integral in the right hand side of (\ref{lem:pvcom3})  under the extra assumption that $\alpha(\epsilon)<1/2$, we prove 
 the following inclusion  
\begin{eqnarray}\label{lem:pvcom5b}
\lefteqn{
\{ \eta\in  {\mathbb{B}}_{n-1}(0,r):\,   \epsilon^2\leq
|\eta|^2+|a\cdot\eta|^2\,,\ 
|\eta|^2+|\gamma(\eta)|^2  <\epsilon^2    \}
}
\\ \nonumber
&&\qquad
\subseteq
\{ \eta\in  {\mathbb{B}}_{n-1}(0,r):\,   \epsilon^2\leq
|\eta|^2+|a\cdot\eta|^2
<
\frac{\epsilon^2}{1-2\alpha(\epsilon)}\}
\end{eqnarray} 
for all $\epsilon\in]0,r[$ such that $\alpha(\epsilon)<1/2$. 
To do so, we note that if  $ \epsilon^2\leq
|\eta|^2+|a\cdot\eta|^2$  and $|\eta|^2+|\gamma(\eta)|^2 <\epsilon^2$, then
\begin{eqnarray*}
\lefteqn{
\epsilon^2 > |\eta|^2+|(\gamma(\eta)-a\cdot\eta)+a\cdot\eta|^2   
}
\\ \nonumber
&&\qquad
\geq |\eta|^2+| \,|\gamma(\eta)-a\cdot\eta|-|a\cdot\eta|\,|^2  
\\ \nonumber
&&\qquad
=|\eta|^2+|a\cdot\eta|^2+| \gamma(\eta)-a\cdot\eta |^2
-2|a\cdot\eta| |\gamma(\eta)-a\cdot\eta|
\\ \nonumber
&&\qquad
\geq |\eta|^2+|a\cdot\eta|^2
-2|a\cdot\eta| |\gamma(\eta)-a\cdot\eta|
\\ \nonumber
&&\qquad
 \geq	 
(|\eta|^2+|a\cdot\eta|^2)
\left(
1-2\frac{
|a\cdot\eta|\alpha(\epsilon)|\eta|
}{|\eta|^2+|a\cdot\eta|^2}
\right)
\geq (|\eta|^2+|a\cdot\eta|^2)(1-2\alpha(\epsilon))
\end{eqnarray*}
 for all $\epsilon\in]0,r[$ such that $\alpha(\epsilon)<1/2$.  Next we turn to estimate the first integral in the right hand side of (\ref{lem:pvcom3}). By the inclusion of (\ref{lem:pvcom4}), we have
 \begin{eqnarray*}
\lefteqn{
\int_{
\{ \eta\in  {\mathbb{B}}_{n-1}(0,r):\,   \epsilon^2\leq |\eta|^2+|\gamma(\eta)|^2\,,
\      |\eta|^2+|a\cdot\eta|^2<\epsilon^2    \}
}|g(\eta)|\,d\eta
}
\\ \nonumber
&&\qquad
\leq c_g
\int_{
\{ \eta\in  {\mathbb{B}}_{n-1}(0,r):\, 
\left(\frac{\epsilon}{1+\alpha(\epsilon)}\right)^2\leq   |\eta|^2+|a\cdot\eta|^2<\epsilon^2    \}
}|\eta|^{-(n-1)}\,d\eta
\end{eqnarray*}
 for all $\epsilon\in ]0, r_g  [$. Now let $R\in O_{n-1}({\mathbb{R}})$ be such that
 \[
a=|a|Re_1\,. 
 \]
Then by changing the variable in the integral with $\eta=R\xi$, we obtain
\begin{eqnarray*}
\lefteqn{
c_g\int_{
\{ \eta\in  {\mathbb{B}}_{n-1}(0,r):\, 
\left(\frac{\epsilon}{1+\alpha(\epsilon)}\right)^2\leq   |\eta|^2+|a\cdot\eta|^2<\epsilon^2    \}
}|\eta|^{-(n-1)
}\,d\eta
}
\\ \nonumber
&&\qquad
=c_g\int_{
\{ \xi\in  {\mathbb{B}}_{n-1}(0,r):\, 
\left(\frac{\epsilon}{1+\alpha(\epsilon)}\right)^2\leq   |\xi|^2+|(|a|e_1)\cdot\xi|^2<\epsilon^2    \}
}|\xi|^{-(n-1)
}\,d\xi
\\ \nonumber
&&\qquad
=c_g\int_{
\{ \xi\in  {\mathbb{B}}_{n-1}(0,r):\, 
\left(\frac{\epsilon}{1+\alpha(\epsilon)}\right)^2\leq   |\xi|^2
\left(1+\xi_1^2\frac{|a|^2}{|\xi|^2}\right)
<\epsilon^2    \}
}|\xi|^{-(n-1)
}\,d\xi
\\ \nonumber
&&\qquad
=c_g\int_{
\left\{ \xi\in  {\mathbb{B}}_{n-1}(0,r):\, 
\left(\frac{\epsilon}{1+\alpha(\epsilon)}\right)\sqrt{1+\xi_1^2\frac{|a|^2}{|\xi|^2}}^{-1}\leq 
|\xi|
<\epsilon\sqrt{1+\xi_1^2\frac{|a|^2}{|\xi|^2}}^{-1}\right\}
}|\xi|^{-(n-1)
}\,d\xi\,.
\end{eqnarray*}
If $n>2$, then an  integration on the spheres implies that the last integral equals
\begin{eqnarray*}
\lefteqn{
c_g\int_{\partial{\mathbb{B}}_{n-1}(0,1)}
\int_{
\left(\frac{\epsilon}{1+\alpha(\epsilon)}\right)\sqrt{1+u_1^2 |a|^2 }^{-1}
}^{
\epsilon\sqrt{1+u_1^2 |a|^2 }^{-1}
}|ru|^{-(n-1)}r^{n-2}\,dr\,d\sigma_u
}
\\ \nonumber
&&\qquad
=c_g\int_{\partial{\mathbb{B}}_{n-1}(0,1)}
\log\left(
\frac{
\epsilon\sqrt{1+u_1^2 |a|^2 }^{-1}
}
{
\left(\frac{\epsilon}{1+\alpha(\epsilon)}\right)\sqrt{1+u_1^2 |a|^2 }^{-1}
}
\right)\,d\sigma_u
\\ \nonumber
&&\qquad
=c_g\int_{\partial{\mathbb{B}}_{n-1}(0,1)}
\log
(1+\alpha(\epsilon))
\,d\sigma_u=c_gs_{n-1}\log
(1+\alpha(\epsilon))
\end{eqnarray*} 
  for all $\epsilon\in ]0, r_g [$. If instead $n=2$, then the last integral above equals
  \[
  2c_g\int_{
 \left(\frac{\epsilon}{1+\alpha(\epsilon)}\right)\sqrt{1+  |a|^2 }^{-1}
}
^{\epsilon\sqrt{1+ |a|^2 }^{-1} }|\xi|^{-(2-1)
}\,d\xi
=c_gs_{1}\log
(1+\alpha(\epsilon))\qquad \epsilon\in ]0, r_g [\,.
  \] 
 Hence,
  \begin{eqnarray}\label{lem:pvcom6}
\lefteqn{
\int_{
\{ \eta\in  {\mathbb{B}}_{n-1}(0,r):\,   \epsilon^2\leq |\eta|^2+|\gamma(\eta)|^2\,,
\      |\eta|^2+|a\cdot\eta|^2<\epsilon^2    \}
}|g(\eta)|\,d\eta
}
\\ \nonumber
&&\qquad\qquad\qquad
 \leq c_gs_{n-1}
\log
(1+\alpha(\epsilon))
 \qquad\forall\epsilon\in ]0, r_g [\,,
\end{eqnarray} 
for all $n\geq 2$, with the understanding that $s_1=2$.  Next we turn to estimate the second integral in the right hand side of (\ref{lem:pvcom3}). 
 By the inclusion of (\ref{lem:pvcom5a}), we have
 \begin{eqnarray*}
\lefteqn{
\int_{
\{ \eta\in  {\mathbb{B}}_{n-1}(0,r):\,   \epsilon^2\leq
 |\eta|^2+|a\cdot\eta|^2\,,
\ 
|\eta|^2+|\gamma(\eta)|^2    <\epsilon^2    \}
}|g(\eta)|\,d\eta
}
\\ \nonumber
&&\quad
\leq
c_g\int_{
\{ \eta\in  {\mathbb{B}}_{n-1}(0,r):\,   \epsilon^2\leq
|\eta|^2+|a\cdot\eta|^2\,,\ 
|\eta|<\epsilon
 \}
}|\eta|^{-(n-1)}\,d\eta
\qquad\forall\epsilon\in ]0, r_g [\,,
\end{eqnarray*} 
 and by the same argument above, we have
 \begin{eqnarray}\label{lem:pvcom7a}
\lefteqn{
\int_{
\{ \eta\in  {\mathbb{B}}_{n-1}(0,r):\,   \epsilon^2\leq
 |\eta|^2+|a\cdot\eta|^2\,,
\ 
|\eta|^2+|\gamma(\eta)|^2   <\epsilon^2    \}
}|g(\eta)|\,d\eta
}
\\ \nonumber
&&\qquad
\leq
c_g \int_{\partial{\mathbb{B}}_{n-1}(0,1)}
\log\left(
\frac{
 \epsilon
}
{
\epsilon\sqrt{1+u_1^2 |a|^2 }^{-1}
}
\right)\,d\sigma_u
\\ \nonumber
&&\qquad
=c_gs_{n-1}\log \sqrt{1+  |a|^2 }
\qquad\forall\epsilon\in ]0, r_g [ 
 \end{eqnarray}
 for $n>2$ and
 \begin{eqnarray}\label{lem:pvcom7aa}
\lefteqn{
 \int_{
\{ \eta\in  {\mathbb{B}}_{n-1}(0,r):\,   \epsilon^2\leq
 |\eta|^2+|a\cdot\eta|^2\,,
\ 
|\eta|^2+|\gamma(\eta)|^2   <\epsilon^2    \}
}|g(\eta)|\,d\eta
}
\\ \nonumber
&&\qquad
=2c_g\int_{\epsilon\sqrt{1+ |a|^2 }^{-1}}^\epsilon\eta^{-1}\,d\eta
\leq c_gs_{1}\log \sqrt{1+  |a|^2 } 
\qquad\forall\epsilon\in ]0, r_g [ 
\end{eqnarray}
  if $n=2$. Next we turn to estimate the second integral in the right hand side of (\ref{lem:pvcom3}) under the extra assumption that $\alpha(\epsilon)<1/2$.  
 By the inclusion of (\ref{lem:pvcom5b}), we have
 \begin{eqnarray*}
\lefteqn{
\int_{
\{ \eta\in  {\mathbb{B}}_{n-1}(0,r):\,   \epsilon^2\leq
 |\eta|^2+|a\cdot\eta|^2\,,
\ 
|\eta|^2+|\gamma(\eta)|^2    <\epsilon^2    \}
}|g(\eta)|\,d\eta
}
\\ \nonumber
&&\quad
\leq
c_g\int_{
\{ \eta\in  {\mathbb{B}}_{n-1}(0,r):\,   \epsilon^2\leq
|\eta|^2+|a\cdot\eta|^2
<
\frac{\epsilon^2}{1-2\alpha(\epsilon)}\}
}|\eta|^{-(n-1)}\,d\eta
\end{eqnarray*} 
 for all $\epsilon\in ]0,r_g[$ such that $\alpha(\epsilon)<1/2$ 
 and by the same argument above, we have
 \begin{eqnarray}\label{lem:pvcom7b}
\lefteqn{
\int_{
\{ \eta\in  {\mathbb{B}}_{n-1}(0,r):\,   \epsilon^2\leq
 |\eta|^2+|a\cdot\eta|^2\,,
\ 
|\eta|^2+|\gamma(\eta)|^2    <\epsilon^2    \}
}|g(\eta)|\,d\eta
}
\\ \nonumber
&&\qquad
\leq
c_g \int_{\partial{\mathbb{B}}_{n-1}(0,1)}
\log\left(
\frac{
 \frac{\epsilon}{\sqrt{1-2\alpha(\epsilon)}} \sqrt{1+u_1^2 |a|^2 }^{-1}
}
{
\epsilon\sqrt{1+u_1^2 |a|^2 }^{-1}
}
\right)\,d\sigma_u
\\ \nonumber
&&\qquad
=c_gs_{n-1}\log\frac{1}{ \sqrt{1-2\alpha(\epsilon)} }
 \end{eqnarray}
for $n>2$ and
\begin{eqnarray}\label{lem:pvcom7bb}
\lefteqn{
\int_{
\{ \eta\in  {\mathbb{B}}_{n-1}(0,r):\,   \epsilon^2\leq
 |\eta|^2+|a\cdot\eta|^2\,,
\ 
|\eta|^2+|\gamma(\eta)|^2    <\epsilon^2    \}
}|g(\eta)|\,d\eta
}
\\ \nonumber
&&\qquad
\leq 2c_g\int_{\epsilon\sqrt{1+ |a|^2 }^{-1}}^{\frac{\epsilon}{\sqrt{1-2\alpha(\epsilon)}} \sqrt{1+  |a|^2 }^{-1}}\eta^{-(2-1)}\,d\eta
=c_gs_{1}\log\frac{1}{ \sqrt{1-2\alpha(\epsilon)} }
\end{eqnarray}
for $n=2$ and  for all $\epsilon\in ]0,r_g[$ such that $\alpha(\epsilon)<1/2$ .
 By inequalities (\ref{lem:pvcom6}) and  (\ref{lem:pvcom7a}), (\ref{lem:pvcom7aa}), we deduce the validity of inequality (\ref{lem:pvcom1b}). By inequalities (\ref{lem:pvcom6}) and  (\ref{lem:pvcom7b}), (\ref{lem:pvcom7bb}), we deduce the validity of inequality (\ref{lem:pvcom1c}).

 By the differentiability of $\gamma $ at $0$, the limiting relation 
(\ref{lem:pvcom2}) holds true  and thus
$\lim_{\epsilon\to0}\alpha(\epsilon)=0$. Hence, there exists $\epsilon_0\in]0,r_g[$ such that 
\[
 \alpha(\epsilon)<1/2\qquad\forall\epsilon\in]0,\epsilon_0[\,.
\]
Then   inequality (\ref{lem:pvcom1c}) and the limiting relation $\lim_{\epsilon\to0}\alpha(\epsilon)=0$ imply
  the validity of the limiting relation 
  \[
\lim_{\epsilon\to 0} \left[
\int_{{\mathbb{B}}_{n-1}(0,r)\setminus A_\gamma(\epsilon)}g(\eta)\,d\eta
- \int_{{\mathbb{B}}_{n-1}(0,r)\setminus A_a(\epsilon)}g(\eta)\,d\eta\right]=0\,.
\]
Hence,  also the last part of statement (iii) holds true and the proof is complete.
\hfill  $\Box$ 

\vspace{\baselineskip}

\noindent
{\bf Acknowledgement:} The author  acknowledges  the support of the Research  
Project GNAMPA-INdAM   $\text{CUP}\_$E53C22001930001 `Operatori differenziali e integrali in geometria spettrale' and the Project funded by the European Union – Next Generation EU under the National Recovery and Resilience Plan (NRRP), Mission 4 Component 2 Investment 1.1 -Call for tender PRIN 2022 No. 104 of February, 2 2022 of Italian Ministry of University and Research; Project code: 2022SENJZ3 (subject area: PE - Physical Sciences and Engineering) "Perturbation problems and asymptotics for elliptic differential equations: variational and potential theoretic methods".

\end{document}